\documentclass[a4paper]{paper}
\usepackage{amsmath,amsthm,amsfonts,mathrsfs}
\usepackage{amssymb}
\usepackage{tikz-cd}
\usepackage{xspace}
\usepackage{graphicx}
\usepackage{paralist}
\usepackage{hyperref} 

\usepackage{multirow}
\usepackage{diagbox}

\usepackage{listings}

\hypersetup{%
    pdfmenubar=true,       
    pdffitwindow=false,     
    pdfstartview={FitH},    
    pdftitle={A Fast and Adaptive Algorithm to Compute the X-ray Transform},            
    colorlinks=true,       
    linkcolor=red,          
    citecolor=red,        
    filecolor=magenta,      
    urlcolor=cyan,           
    pdfborder = {0,0,0}
}

\usepackage{xspace}

\usepackage{enumerate}
\usepackage{subfig}
\usepackage{booktabs}

\usepackage{acro}
\usepackage{cleveref}
\usepackage{algorithm}
\usepackage[noend]{algpseudocode}
\makeatletter
\def\BState{\State\hskip-\ALG@thistlm}
\makeatother

\Crefname{equation}{}{}
\crefname{equation}{}{}

\usepackage{todonotes}

\usepackage{specialmacros}
\usepackage{paperacronyms}

\begin{document}

\title{A Fast and Adaptive Algorithm to Compute the X-ray Transform}
\author{
Chong Chen\thanks{LSEC, ICMSEC, Academy of Mathematics and Systems Science, Chinese Academy of Sciences, Beijing 100190, China (chench@lsec.cc.ac.cn).}, Runqian Wang\thanks{Princeton International School of Mathematics and Science, Princeton, NJ 08540, USA (ray.wang@prismsus.org).}, 
Chandrajit Bajaj\thanks{Department of Computer Science, The Oden Institute for Computational Engineering and Sciences, University of Texas at Austin, Austin, TX 78712, USA (bajaj@cs.utexas.edu).} and 
Ozan \"Oktem\thanks{Department of Mathematics, KTH - Royal Institute of Technology, 100 44 Stockholm, Sweden (ozan@kth.se).} 
}

\maketitle

\begin{abstract}
We propose a new algorithm to compute the X-ray transform of an image represented by unit (pixel/voxel) basis functions. The fundamental issue is equivalently calculating the intersection lengths of the ray with associated units. For any given ray, we first derive the sufficient and necessary condition for non-vanishing intersectability. By this condition, we then distinguish the units that produce valid intersections with the ray. Only for those units rather than all the individuals, we  calculate the intersection lengths by the obtained analytic formula. The proposed algorithm is adapted to 2D/3D parallel beam and 2D fan beam. Particularly, we derive the transformation formulas and generalize the algorithm to 3D circular and helical cone beams. Moreover, we discuss the intrinsic ambiguities of the problem itself, and present a solution. The algorithm not only possesses the adaptability with regard to the center position, scale and size of the image, but also is suited to parallelize with optimality. The comparison study demonstrates the proposed algorithm is fast, more complete, and is more flexible with respect to different scanning geometries and different basis functions. Finally, we validate the correctness of the algorithm by the aforementioned scanning geometries. 
\end{abstract}

\begin{keywords}
X-ray transform, intersection length, projection matrix, sufficient and necessary condition, non-vanishing intersectability, ambiguity and adaptability, tomographic image reconstruction 
\end{keywords}

\section{Introduction}
\label{sec:Introduction}

The tomography is quite critical in clinical diagnosis, such as \ac{X-ray CT}, \ac{ECT}, etc. \cite{hsieh2009computed,cherry2012physics}, and in structure biology as well, such as \ac{Cryo-EM}, \ac{ET}, etc. \cite{frank2006electron,frank2006three}. The X-ray transform models a forward projection operator of image formation for the above imaging modalities, which has been widely used for tomographic image reconstruction \cite{De1983,natterer2001mathematics,frank2006three}.  To reduce the radiation or conduct fast scanning, the low-dose or sparse-view sampling is often required. However, the conventional analytic methods (filtered/weighted backprojection) for image reconstruction does not work well for this kind of data. Recently, the advanced methods based on regularization, compressed sensing and/or deep learning have been gained extensive study, due to their ability to allow reducing dose or scanning views while maintaining or improving reconstructed image quality  (for instance, see \cite{delaney1998globally,JiangMingb,SiKaPa06,thibault2007three,li2011single,hsieh2013recent,dong2013x,sidky2012convex,ramani2012splitting,
xu2014blended,chen_xu_2016,adler2017solving,ArMaOkSc2019,ChGrOz19}). 

It is well-known that the most computationally intensive and tough components for these approaches locates in the computations of X-ray transform and its adjoint \cite{SiKaPa06,ramani2012splitting,chen_xu_2016,ArMaOkSc2019}.  Mathematically, the computation of the adjoint can be converted into calculating X-ray transform. The aim of this work is to investigate the algorithm for computing this transform. 

The X-ray transform is equivalently called \ac{1D} Radon transform \cite{helgason1999radon}. More precisely, the X-ray transform $\ForwardOp$ is given as follows: if $\theta \in \mathbb{S}^{d-1}$ for $d = 2$ or $3$, and $\boldsymbol{x} \in \theta^{\bot}$, then 
\begin{equation}\label{eq:Xray_transform}
\ForwardOp\!f(\theta, \boldsymbol{x}) = \int_{-\infty}^{+\infty} f(\boldsymbol{x} + t\theta)\dint t.
\end{equation}
The formula above is actually the integral of function $f$ over the straight line (ray) through point $\boldsymbol{x}$ with direction $\theta$, which models the forward projection for biomedical imaging modalities, including \ac{X-ray CT}, \ac{Cryo-EM} and \ac{ET} in \cite{natterer2001mathematics,frank2006three} and \cite[Chapter 8]{frank2006electron}. 

Moreover, for \ac{ECT}, the forward projection is mathematically represented by the generalized Radon transform (see \cite[Section II.6--II.7]{natterer2001mathematics} and \cite[Section 4.2]{kak2001principles})
\begin{equation}\label{eq:generalizedRadontransform}
\mathcal{R}_{\Phi}f(\theta, \boldsymbol{x}) = \int_{-\infty}^{+\infty}\Phi(\theta, \boldsymbol{x}, t)f(\boldsymbol{x} + t\theta)\dint t.
\end{equation}
Specifically, the function $\Phi := \Phi(\theta, \boldsymbol{x})$ is independent of the integral variable $t$ for \ac{PET}. Consequently,   \cref{eq:generalizedRadontransform} is reduced to \cref{eq:Xray_transform} readily. In terms of \ac{SPECT}, the function $\Phi := \Phi(\theta, \boldsymbol{x}, t)$ is an exponential attenuation function, and then \cref{eq:generalizedRadontransform} boils down to the attenuated Radon transform, which is a special case of the generalized Radon transform. For both cases, numerically, the function $\Phi$ is given explicitly, and can be used to perform the attenuation correction for image reconstruction. Therefore, the computation of the X-ray transform in \cref{eq:Xray_transform} can be also served as the algorithm for computing the forward projection in \ac{ECT}. 

Note that the function $f$ is the image to be reconstructed, which is often represented by the unit (pixel/voxel) basis functions $\{b_i\}$ as 
\begin{equation}\label{eq:f_represent}
f(\boldsymbol{r}) = \sum_{i}f_ib_i(\boldsymbol{r}),
\end{equation}
where $f_i$ is the given gray value of the $i$-th unit, and 
\[
b_i(\boldsymbol{r}) =  \begin{cases}
  1, &  \boldsymbol{r} \in \Omega_{i}, \\
  0, & \text{otherwise}. 
 \end{cases}
\]
Here $\Omega_{i} \subset \Real^d$ is the $i$-th unit of the image. By \cref{eq:Xray_transform} and \cref{eq:f_represent}, the X-ray transform of the image function can be written as 
\begin{equation}\label{eq:Xray_transform_image}
\ForwardOp\!f(\theta, \boldsymbol{x}) = \sum_{i}f_i\int_{-\infty}^{+\infty} b_i(\boldsymbol{x} + t\theta)\dint t, 
\end{equation}
which is to compute the summation of the intersection lengths of the ray with each unit weighted by its gray value. Hence the fundamental problem of \cref{eq:Xray_transform_image} is to  calculate the line integral of the unit basis functions, which is equivalently the intersection length of the ray with the associated unit. The focus of this article is on studying this problem. 

The intersection lengths of a certain ray with all the units construct a row of the forward projection matrix, where its entry is the associated intersection length. Actually, the computation of the X-ray transform of all rays (forward projection) is also equivalent to compute the projection matrix multiplying a vectorized image, and the computation of their adjoints (backprojection) means to compute the transpose of the projection matrix multiplying a corresponding vector. The projection matrix is of particularly interest to develop and validate the related reconstruction algorithms, which can be stored in sparse pattern if permitted to avoid its repeat calculations during algorithm implementation. 

There are several algorithms developed to compute the forward projection for tomographic imaging, including pixel-driven methods, distance-driven methods and ray-driven methods (see \cite{siddon_1985,ZhGoHe94,JaLe98,HaLi99,ZhRa03,DeBa04,LoFeBa10,gao_2012,LiYu17}). The proposed method of this paper is attributed to the ray-driven method. Such kind of method is often used to model the discrete forward projection \cite{SiKaPa06,ramani2012splitting,sidky2012convex,chen_xu_2016,BuDiSi17}. There are also many software packages that implement the forward projection by ray-driven method, for instance, ASTRA \cite{VaPa16}, Matlab Image Processing Toolbox, etc. A classical accurate algorithm was proposed by Siddon in \cite{siddon_1985}, which requires to compute the intersection points of the ray with all of the grid lines/planes in \ac{2D}/\ac{3D} circumstances, and then sorts all those intersection points. Some speedup versions of Siddon's algorithm were developed in \cite{JaLe98,HaLi99,gao_2012}. Many of them calculate the X-ray transform and its adjoint following Siddon's algorithm, and/or support a few of commonly scanning geometries.

\paragraph{Contributions.} 
In this work, the main contribution is that we propose a fast, accurate, adaptive and parallelizable algorithm to compute the X-ray transform of an image represented by unit (pixel/voxel) basis functions. We first propose the algorithm for the most basic scanning geometry (\ie \ac{2D} parallel beam), and then extend to \ac{2D} fan beam by the associated transformation formula. We further generalize the algorithm to \ac{3D} parallel beam, and derive the transformation formulas and adapt the algorithm to \ac{3D} equiangular and equispaced circular/helical cone beam. The proposed algorithm is not just limited to the above scanning geometries, which is applicable to any beam that its parameters can be transformed into those of the corresponding parallel beam. Additionally, this algorithm is also applicable to the case of the image represented by polygonal/polyhedral basis functions.  So the proposed algorithm is more flexible, \eg one can easily work with different geometries and different image basis functions. 

More importantly, we give the sufficient and necessary condition for non-vanishing intersectability of the ray and unit, which can be used to distinguish out the units that produce valid intersections with the given ray. And then we just calculate the intersection lengths by the obtained analytic formula for those units rather than all the individuals. Based on the results above, the algorithm becomes quite easy to be understood and implemented, and its computational cost is significantly reduced, since we only need to deal with the units that have non-vanishing intersections with the given ray, and the number of non-vanishingly intersectable units is quite few compared to the total number of individuals. The computational cost is $O(N)$ for any given ray, and $O(NM)$ for all rays, to both \ac{2D} and \ac{3D} circumstances. Here the $N$ represents the size of the image being reconstructed along one axis, and $M$ denotes the number of rays. Hence, the computational complexity of this algorithm is optimal. 

We further discuss the intrinsic ambiguities of the problem itself that perhaps happen, and give a solution in the algorithm. Moreover, the algorithm possesses the adaptability with regard to the center position, scale and size of the image, and the more general scanning geometries. Since the proposed algorithm implements the calculation of intersection length ray by ray and unit by unit, and all of the calculations are independent of each other, it is quite suited to parallelize and the computational complexity per parallel thread achieves $O(1)$. The projection matrix can be sparsely stored and output if needed, and the adjoint of X-ray transform can be also computed by the algorithm. Hence, the algorithm is more complete to use, and can be customized freely according to the requirements of the users, and more scanning geometries can be easily added into the framework based on the proposed algorithm.

\paragraph{Outline.} 
The algorithms for various \ac{2D} and \ac{3D} scanning geometries are proposed in \cref{sec:algorithm_2D} and \cref{sec:algorithm_3D}, respectively. \Cref{sec:Discussions} points out the intrinsic ambiguities of the problem itself, and presents a solution accordingly, and also includes the discussions on the adaptability, computational cost and parallelization of the proposed algorithm, and the comparison with existing approaches as well. The validations are performed in \cref{sec:numerical_tests}. Finally, \cref{sec:Conclusions} concludes the paper.

\section{Algorithm for \ac{2D} scanning geometries}
\label{sec:algorithm_2D}

The purpose of this section is to develop the algorithm for \ac{2D} scanning geometries. To begin with, we need to introduce several requisite preliminaries. 

\subsection{Preliminaries}
\label{sec:Preliminaries}

Here we present some preliminaries, including the commonly used \ac{2D} scanning geometries, \ac{2D} imaging coordinate system, and pixel indexes. 

\begin{figure}[htpb]
    \centerline {
    \begin{tabular}{cc}
    \includegraphics[width=0.48\textwidth]{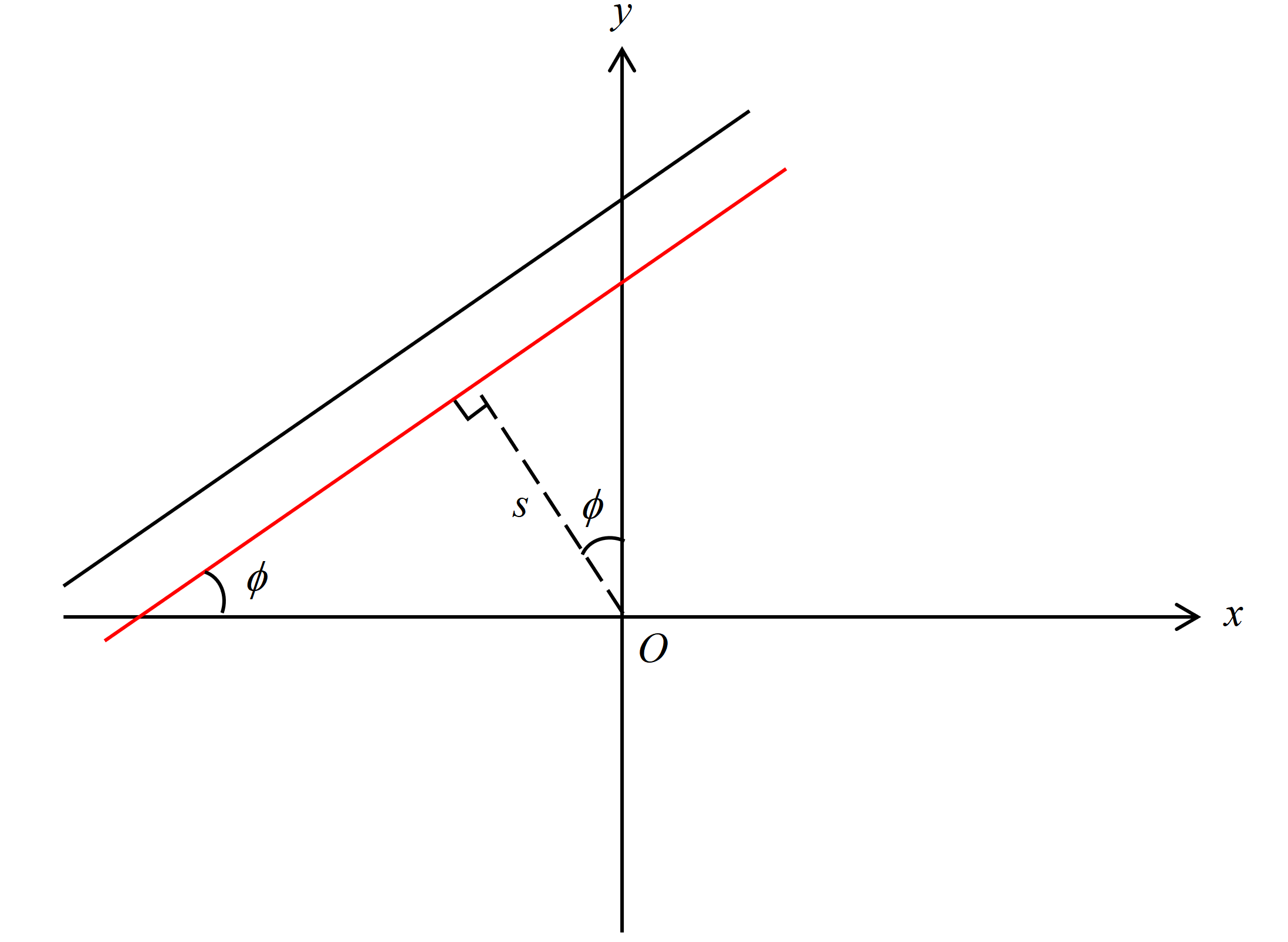}&
    \includegraphics[width=0.48\textwidth]{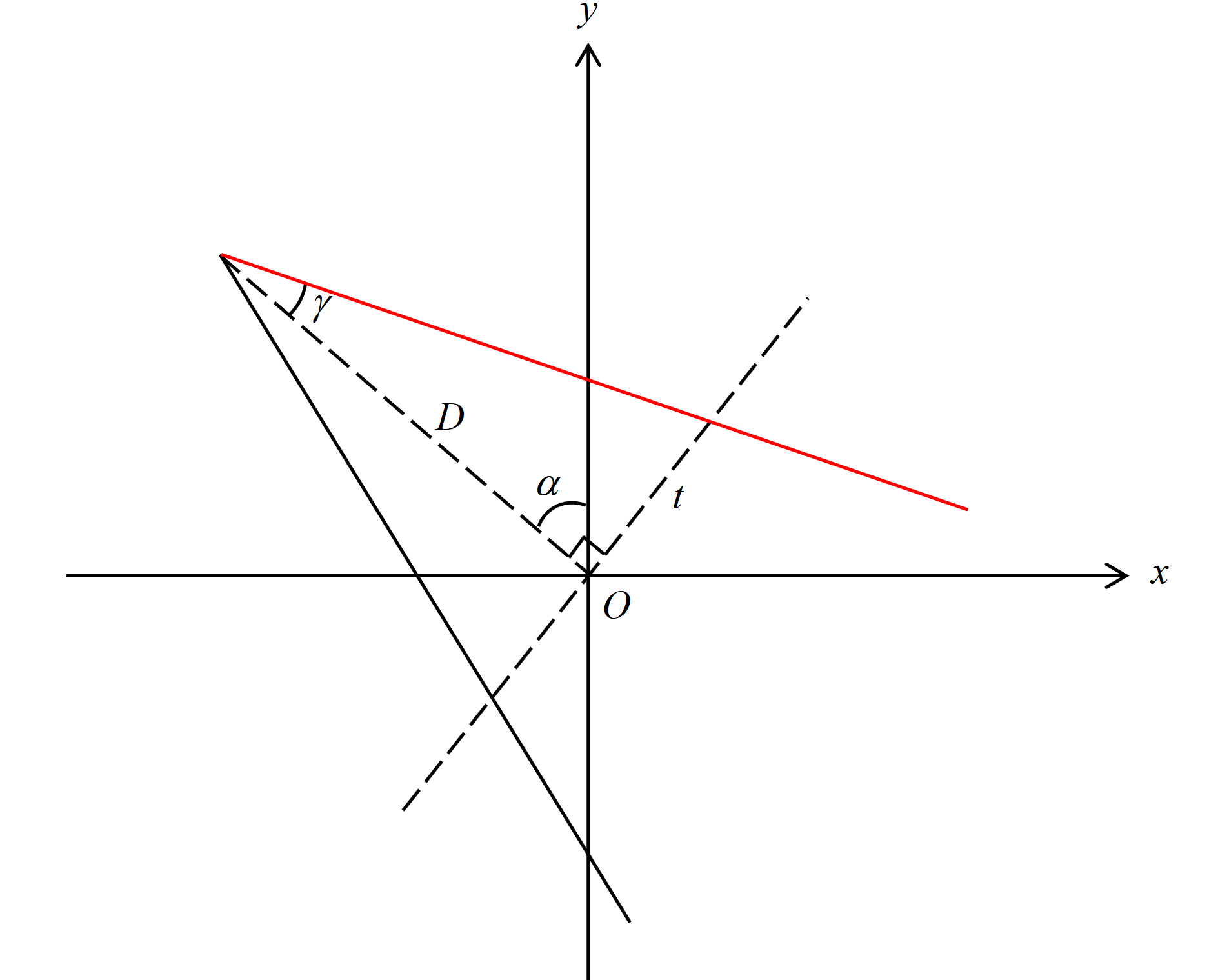}\\
    (a) &  (b) 
    \end{tabular}
    }
    \caption{(a) \ac{2D} parallel beam; (b) Fan beam}
    \label{fig:parallelbeam_fanbeam}
\end{figure}

\paragraph{\ac{2D} scanning geometries.} 
For \ac{2D} case, there are normally two scanning geometries at which the rays are distributed, namely, parallel beam and fan beam (including equiangular and equispaced) \cite{kak2001principles,hsieh2009computed}. Parallel beam, as its name suggested, stands for a group of parallel X-rays penetrating through the detected object, as illustrated in  \cref{fig:parallelbeam_fanbeam} (a). Fan beam, on the other hand, is a set of beams emitted from one source point, and gets detected from a series of detectors, as shown in \cref{fig:parallelbeam_fanbeam} (b). When the detectors are placed such that the angles between any two consecutive rays are equal, it is called equiangular fan beam; When the detectors are aligned with equal space between any two adjacent detectors, it is named equispaced fan beam. 

\begin{figure}[htpb]
    \centerline {
    \begin{tabular}{cc}
    \includegraphics[width=0.45\textwidth]{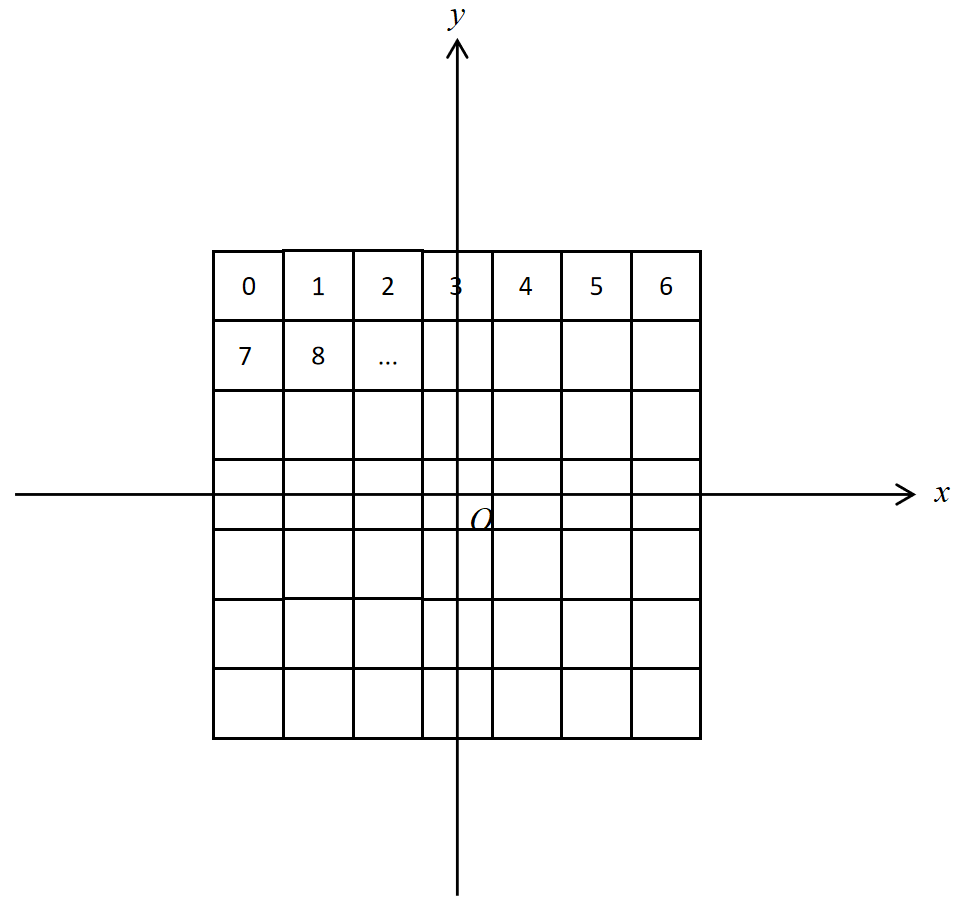}&
    \includegraphics[width=0.45\textwidth]{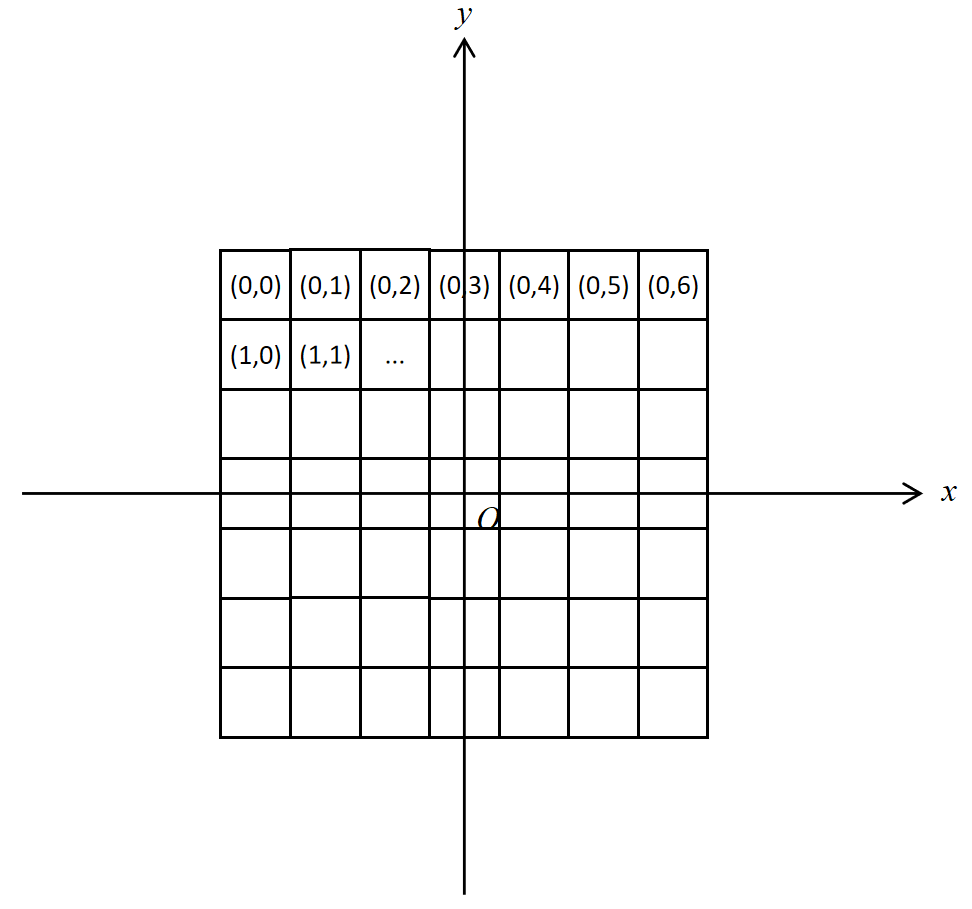}\\
    (a) &  (b)
    \end{tabular}
    }
    \caption{(a) Pixel index in \ac{1D} form; (b) Pixel index in \ac{2D} form}
    \label{fig:2Dindex_2Dindex2Dform}
\end{figure}

\paragraph{\ac{2D} imaging coordinate system.} 
Let $(x, y)$ be the coordinate system. For a image to be reconstructed, we assume that the parameters $L_x$, $L_y$, and $N_x$, $N_y$ are given, where $(L_x, L_y)$ and $(N_x, N_y)$ specify the side lengths of the domain and the size of the image along $x$- and $y$-axis, respectively. Let $(d_x, d_y) := (L_x/N_x, L_y/N_y)$ be the side lengths or scales of the pixel along $x$- and $y$-axis. Without loss of generality, we assume that $L_x = L_y$, $N_x = N_y$, accordingly, $d_x = d_y$, and the center of the image domain is at the origin $O$ of the coordinate system as displayed in \cref{fig:2Dindex_2Dindex2Dform}. For simplicity, we further assume that the scale $d_x = d_y = 1$. If the scale is not unity, the real value of X-ray transform just equals to the scale multiplying that value for the case with unity scale.

\paragraph{Pixel indexes.} 
We define two different indexes for the pixels of the image. The one is given in \ac{1D} form as $I = 0, 1, \ldots, N_xN_y - 1$, which is shown as the example with $N_x = N_y = 7$ in  \cref{fig:2Dindex_2Dindex2Dform} (a). The other one is presented in \ac{2D} form as $(j, i)$ with $j = 0, 1, \ldots, N_y - 1$, $i = 0, 1, \ldots, N_x - 1$, which is shown as the example also with $N_x = N_y = 7$ in \cref{fig:2Dindex_2Dindex2Dform} (b). The relation between them can be expressed as $I = jN_y + i$, which can be easily validated by the examples in \cref{fig:2Dindex_2Dindex2Dform}. They are able to convert into each other. 


Here the first and the second indexes are used to indicate the pixel positions in the associated vectorized image and the original image, respectively. They are independent of the coordinate system. The first one is used in sparse storage of projection matrix in numerical implementation of this paper. Note that the second one is different from the coordinate position of the pixel. Later in computation, when referring to any location of the pixel/voxel, we will use its indexes rather than coordinate position. The coordinate positions are only used in the implementation of the algorithm. 

In this section, we will use pixel $(j, i)$ to indicate the pixel with index $(j, i)$. The valid $j$ (or $i$) is for itself satisfying $0 \leq j \leq N_y -1$ (or $0 \leq i \leq N_x -1$). 

\subsection{Derivation of intersection length}
\label{sec:derivation_intersection_2D}

Here we will give the mathematical derivations of intersection length for various \ac{2D} scanning geometries. 

\subsubsection{\ac{2D} parallel beam}
\label{sec:2D_parallel_beam}

To calculate the intersection length, the position of the ray and the pixel need to be determined. Assuming that the parameters $s \in \Real^1$ and $\phi \in [0, \pi)$ are given, where the $(s, \phi)$ determines a ray of the parallel beam in \ac{2D} plane as shown in \cref{fig:parallelbeam_fanbeam} (a). More precisely, the $s$ is the signed distance between the ray and the origin, and the $\phi$ is the angle between the ray and the positive $x$-axis. 

Given ray $(s, \phi)$, the resulting unit direction and normal of the ray are  
\begin{equation}\label{eq:direction_normal}
\theta = (\cos{\phi}, \sin{\phi}) \quad \text{and} \quad \theta^{\bot} = (-\sin{\phi}, \cos{\phi}), 
\end{equation}
respectively. 

As assumed in the above, the scale $d_x = d_y = 1$. We then give the transformation relation between the \ac{2D}-form index of the pixel and the coordinate position of its center. For any pixel with index $(j, i)$, its center is located at the coordinate position 
\[
\biggl(i - \frac{N_x - 1}{2}, ~ \frac{N_y - 1}{2} - j\biggr).
\]
Hence the associated pixel basis function is defined as 
\begin{equation}\label{eq:pixel_basis_function}
    f_{ji}(x, y) = \begin{cases}
  1, & (x, y) \in \Omega_{ji}, \\
  0, & \text{otherwise}, 
 \end{cases}
\end{equation}
where $\Omega_{ji} := \bigl[i - \frac{N_x}{2}, i - \frac{N_x}{2} + 1\bigr] \times \bigl[\frac{N_y}{2} - j - 1, \frac{N_y}{2} - j\bigr]$ denotes the support of pixel $(j, i)$.

The intersection length $l_{ji}(s, \phi)$ between the ray $(s, \phi)$ and pixel $(j, i)$ can be expressed as the X-ray transform of the associated pixel basis function, namely, 
\begin{equation}\label{eq:intersection_length}
    l_{ji}(s, \phi) = \int_{-\infty}^{\infty}f_{ji}(t\theta + s\theta^{\bot})\dint t.
\end{equation}
Using \cref{eq:direction_normal} and \cref{eq:pixel_basis_function}, the \cref{eq:intersection_length} implies that if and only if the following condition is satisfied  
\begin{equation}\label{eq:condition_nonzero}
\begin{cases}
  i - \frac{N_x}{2} \leq t\cos\phi - s\sin\phi \leq i - \frac{N_x}{2} + 1,  \\
  \frac{N_y}{2} - j - 1 \leq t\sin\phi + s\cos\phi \leq \frac{N_y}{2} - j, 
 \end{cases}
\end{equation}
then the length of intersection might be nonzero. Note that the variable $t$ of integration in \cref{eq:intersection_length} is along the ray. Therefore, the range of $t$ that satisfies the inequalities in \cref{eq:condition_nonzero} would be the length of intersection. 

As observed, the ray $(s, \phi + 2n\pi)$ for $n\in \mathbb{Z}$ is equivalent to the ray $(s, \phi)$, and the ray $(s, \phi + \pi)$ is just the same as the ray $(-s, \phi)$. Hence any ray can be determined by the ray $(s, \phi) \in \Real^1 \times [0, \pi)$. In what follows we list the cases by the choice of $\phi$. 

\paragraph{Case 1: $\phi \in (0, \pi/2)$.}
For simplicity, let 
\begin{equation}\label{eq:C_xC_y}
C_x = \frac{i - N_x/2 + s\sin\phi}{\cos\phi}, \quad C_y = \frac{N_y/2 - j - s\cos\phi}{\sin\phi}.
\end{equation}
Using simple calculations, by \cref{eq:C_xC_y}, the \cref{eq:condition_nonzero} can be translated into 
\begin{equation}\label{eq:condition_nonzero_t}
\begin{cases}
  C_x \leq t  \leq C_x + \frac{1}{\cos\phi},  \\
  C_y - \frac{1}{\sin\phi} \leq t \leq C_y. 
 \end{cases}
\end{equation}
As we observed, if the two intervals in  \cref{eq:condition_nonzero_t} are disjointed or their overlap is just a point, then the length of the intersection is definitely vanishing. To exclude this possibility, if and only if 
\begin{equation}\label{eq:conditions_2D_1}
C_{low} < C_{up}, 
\end{equation}
where 
\[
C_{low} = \max(C_x, C_y - 1/\sin\phi), \quad C_{up} = \min(C_x + 1/\cos\phi, C_y).
\] 
Evidently, \cref{eq:conditions_2D_1} is equivalent to 
\begin{equation}\label{eq:conditions_2D_2}
\begin{cases}
 C_x < C_y, \\
 C_y - \frac{1}{\sin\phi} < C_x + \frac{1}{\cos\phi}. 
 \end{cases}
\end{equation}
The \cref{eq:conditions_2D_1} or \cref{eq:conditions_2D_2} is just the \textit{sufficient and necessary condition} for non-vanishing intersectability of the given ray and pixel. More precisely, the \textit{non-vanishing intersectability} means the length of the intersection being non-vanishing. 

More specifically, by \cref{eq:C_xC_y}, if the $i$ is given, the condition \cref{eq:conditions_2D_2} is equivalent to the following inequality 
\begin{equation}\label{eq:restriction_j}
C_1  - \tan\phi - 1 < j < C_1, 
\end{equation}
where
\[        
C_1 = \frac{N_y}{2} - s\cos\phi - C_x\sin\phi.
\]
Similarly, provided that the $j$ is given, the condition \cref{eq:conditions_2D_2} can be rewritten as 
\begin{equation}\label{eq:restriction_i}
C_2  - \cot\phi - 1 < i < C_2, 
\end{equation}
where
\[
    C_2 = \frac{N_x}{2} - s\sin\phi + C_y\cos\phi.
\]
Hence, for any given $i$ or $j$, we can calculate the range of valid $j$ by \cref{eq:restriction_j}, or that of valid $i$ by \cref{eq:restriction_i}, to obtain those pixels intersecting with the given ray non-vanishingly. For those pixels, the condition \cref{eq:conditions_2D_1} is surely satisfied. By merging the intervals in \cref{eq:condition_nonzero_t}, we get the intersection as 
\begin{equation}\label{eq:range_t}
   C_{low} \leq t \leq C_{up}. 
\end{equation}
Otherwise, the intersection is empty. Hence, the \textit{analytic formula} of the non-vanishing intersection length should be 
\begin{equation}\label{eq:range_length}
    l_{ji}(s, \phi) = C_{up} - C_{low}.
\end{equation}
One alternative analytic formula was derived for the Radon transform over a fixed unit square in \cite[Chapter 2]{De1983}, which is given under the condition of the different ranges for the distance $s$. Here the analytic formula is presented without that restriction. Compared to that counterpart, the derivation of \cref{eq:range_length} is simpler and more general, and much easier to be generalized into higher-dimensional scenarios (see \cref{sec:algorithm_3D}). 

In what follows we explicitly give a sketch of the algorithm for $\phi \in (0, \pi/2)$. 
\begin{algorithm}
\caption{A sketch of the algorithm for $\phi \in (0, \pi/2)$.}
\label{algo:algorithm_2D}
\begin{algorithmic}[1]
\State Given $(s, \phi)$ to determine a ray.
\State \emph{Loop}: For $0 \leq i \leq N_x - 1$ 
\State \quad\quad\quad\quad\quad Obtain the valid $j$ by \cref{eq:restriction_j}. 
\State \quad\quad\quad\quad\quad Calculate the intersection length by \cref{eq:range_length} for pixel ($j, i$). 
\State \emph{End loop}.
\State \textbf{Output} the intersection lengths and associated \ac{1D}-form indexes.
\end{algorithmic}
\end{algorithm}

\paragraph{Case 2: $\phi = 0$.} The \cref{eq:condition_nonzero} becomes 
\begin{equation}\label{eq:condition_nonzero_case3}
\begin{cases}
  i - \frac{N_x}{2} \leq t \leq i - \frac{N_x}{2} + 1,  \\
  \frac{N_y}{2} - j - 1 \leq s \leq \frac{N_y}{2} - j.
 \end{cases}
\end{equation}
Obviously, for any given $i$, if the valid $j$ is satisfying $\frac{N_y}{2} - s - 1 \leq j \leq \frac{N_y}{2} - s$, 
the intersection length is unity, otherwise, it is zero. 

\begin{remark}\label{rem:cases_2D}
For the case $\phi \in (\pi/2, \pi)$, the corresponding sufficient and necessary condition for non-vanishing intersectability, and analytic formula of the intersection length can be easily obtained by making quite slight changes to the derivation in the above \textbf{case 1}. 

Similarly, the intersection length can be immediately calculated for the case $\phi = \pi/2$ by the method of \textbf{case 2}. Note that the rays in such cases are parallel to a certain axis, which results in quite simple calculations, but the ambiguity perhaps happen (see the first example in \cref{sec:ambiguity_solution}). 
\end{remark}

\subsubsection{Fan beam}
\label{sec:fan_beam}

Apart from the \ac{2D} parallel beam, as mentioned previously, the fan beam is often used in \ac{2D} circumstance. Therefore, the algorithm should be able to adapt to this situation. 

\paragraph{Equiangular fan beam.} The second pattern is equiangular fan beam as shown in \cref{fig:parallelbeam_fanbeam} (b). Essentially, any ray of the equiangular fan beam can be determined by the parameters $(D, \alpha, \gamma)$,  where $D > 0$ indicates the distance of source for X-rays from the origin $O$, $\alpha \in [0, 2\pi)$ stands for the angle between the line from the origin to source and the positive $y$-axis, and $\gamma \in [-\gamma_m, \gamma_m]$ specifies the signed angle between the line from the source to the origin and each ray. The radian $\gamma_m \in (0, \pi/2)$ determines the size of the view field for fan beam.   

Then the aim is to find the corresponding parameters $(s, \phi)$ in \ac{2D} parallel beam for each ray when given $(D, \alpha, \gamma)$, as indicated in \cref{fig:parallelbeam_fanbeam} (b). In this case, it is quite straightforward:
\begin{equation}\label{eq:2Dparalleltofan}
s = D\sin\gamma, \quad \phi = \gamma + \alpha - \frac{\pi}{2}.  
\end{equation}
The above transformation can be also referred to \cite{kak2001principles}. After transformation, the associated X-ray transform can be readily computed by the algorithm developed in \cref{sec:2D_parallel_beam}. 
    
\paragraph{Equispaced fan beam.} As depicted in \cref{fig:parallelbeam_fanbeam} (b), the equispaced fan beam is quite similar with the equiangular one. But the one given parameter is about the signed distance along the line corresponding to the detector bank \cite{kak2001principles}. In other words, the given parameter is $t \in [-t_m, t_m]$ compared to the $\gamma$ in the previous case. Here the $t_m > 0$ determines the size of the view field for fan beam. Thus, for each ray with the given parameters $(D, \alpha, t)$, in what follows we only need to perform several minor changes to the derivation above 
\begin{equation}\label{eq:2Dparalleltofan_2}
     s = \frac{Dt}{\sqrt{D^2 + t^2}}, \quad \phi = \arctan\bigl(\frac{t}{D}\bigr) +\alpha-\frac{\pi}{2}.  
\end{equation}
Similarly, the associated X-ray transform in equispaced fan beam can be computed by the algorithm developed in \cref{sec:2D_parallel_beam}.

\section{Algorithm for \ac{3D} scanning geometries}
\label{sec:algorithm_3D} 

Here the aim is to further develop the algorithm for \ac{3D} scanning geometries. To proceed, we first introduce some preliminaries. 

\subsection{Preliminaries}
\label{sec:Preliminaries_3D}

The purpose of this part is to state some requisite preliminaries for \ac{3D} case, including the often used \ac{3D} scanning geometries, \ac{3D} imaging coordinate system, voxel indexes, and Eulerian angles. 

\begin{figure}[htbp]
    \centerline {
    \begin{tabular}{ccc}
    \includegraphics[width=0.32\textwidth]{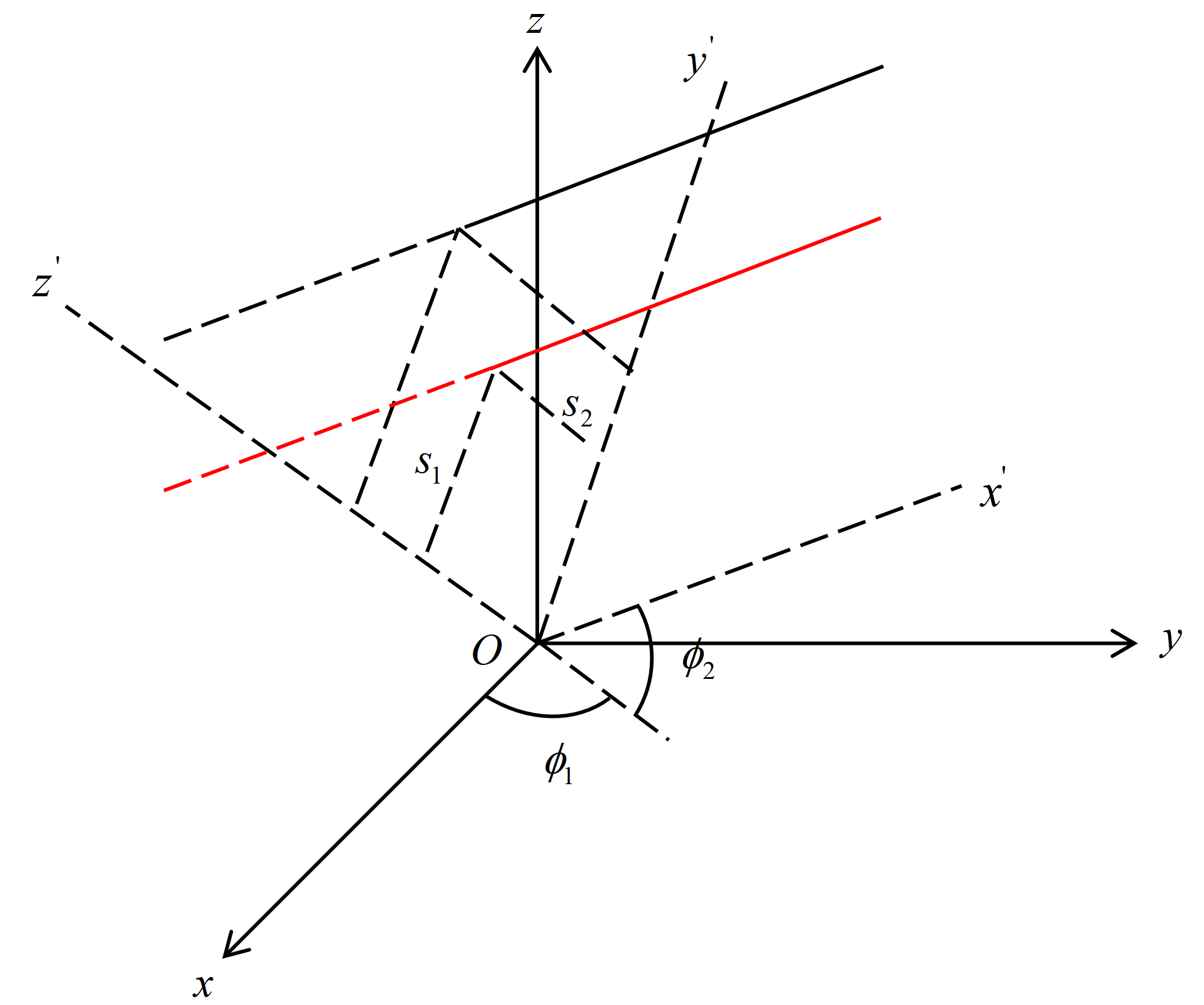}&
    \includegraphics[width=0.32\textwidth]{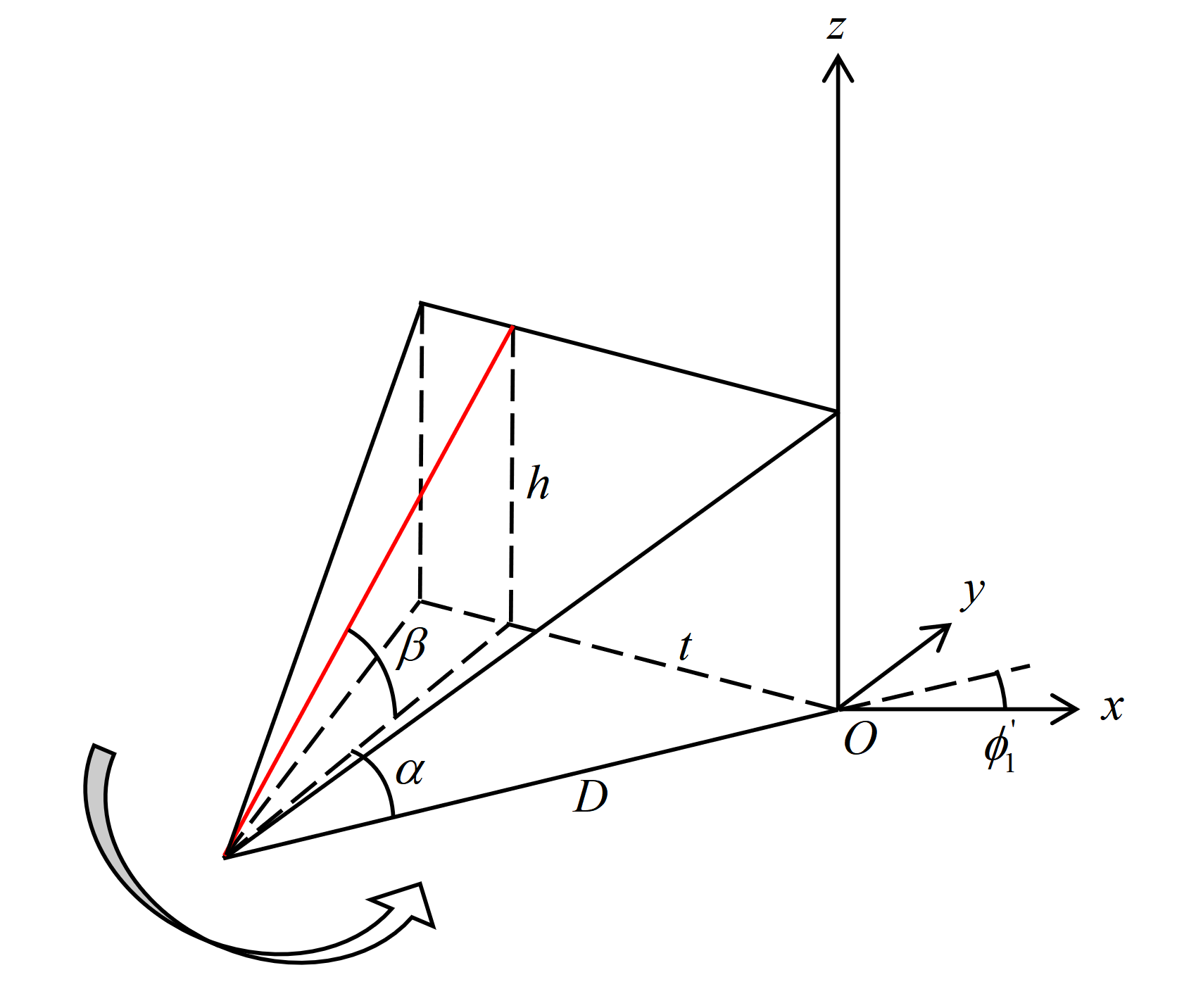}&
    \includegraphics[width=0.32\textwidth]{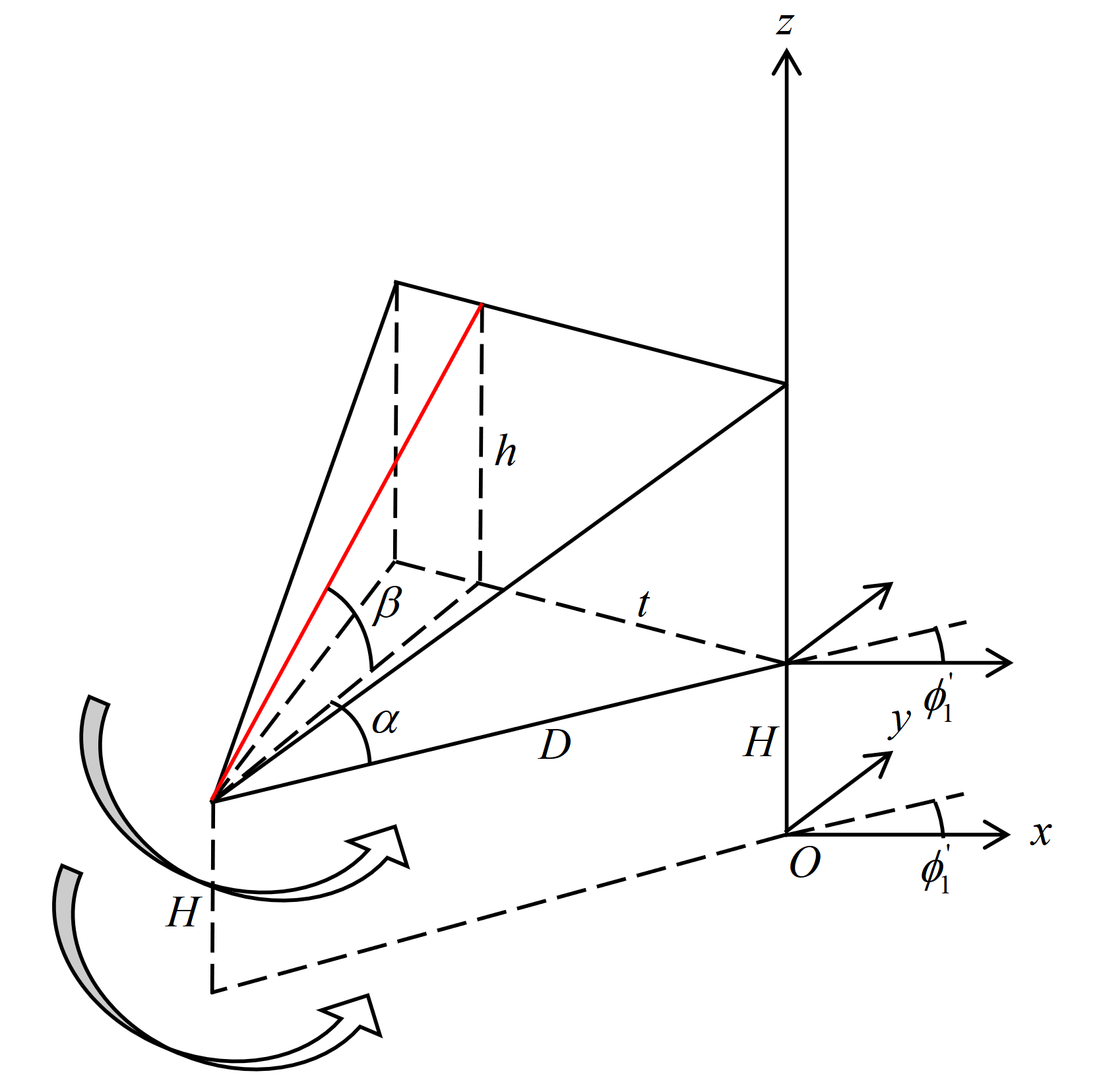}\\
    (a) & (b) & (c)
    \end{tabular}
    }
    \caption{(a) \ac{3D} parallel beam; (b) Circular cone beam; (c) Helical cone beam}
    \label{fig:3D_parallel_beam_circular_cone_beam_Helical_cone_beam}
\end{figure}

\paragraph{\ac{3D} scanning geometries.} 
In \ac{3D} circumstance, there are three patterns commonly used, including parallel beam, circular cone beam and helical cone beam \cite{kak2001principles,frank2006three,hsieh2009computed}. The \ac{3D} Parallel beam is almost the same as in \ac{2D} situation, where all beams from one view are parallel to each other, as illustrated in \cref{fig:3D_parallel_beam_circular_cone_beam_Helical_cone_beam} (a). 
Apart from the applications in medical imaging, such beam has been widely applied in biological imaging (e.g., Cryo-EM and ET) \cite{frank2006three,frank2006electron}.  Circular cone beam can be seen as the extension of the fan beam, where the X-ray source is performing uniformly circular motion around the rotating axis, as depicted in \cref{fig:3D_parallel_beam_circular_cone_beam_Helical_cone_beam} (b). This beam naturally includes the equiangular and equispaced patterns. 

Helical cone beam is often used when scanning through a long object, where the X-ray source relatively moves with constant speed along one axis while performing uniformly circular motion on the perpendicular plane \cite{hsieh2009computed}, as shown in \cref{fig:3D_parallel_beam_circular_cone_beam_Helical_cone_beam} (c). 

\begin{figure}[htbp]
    \centerline {
    \begin{tabular}{cc}
    \includegraphics[width=0.48\textwidth]{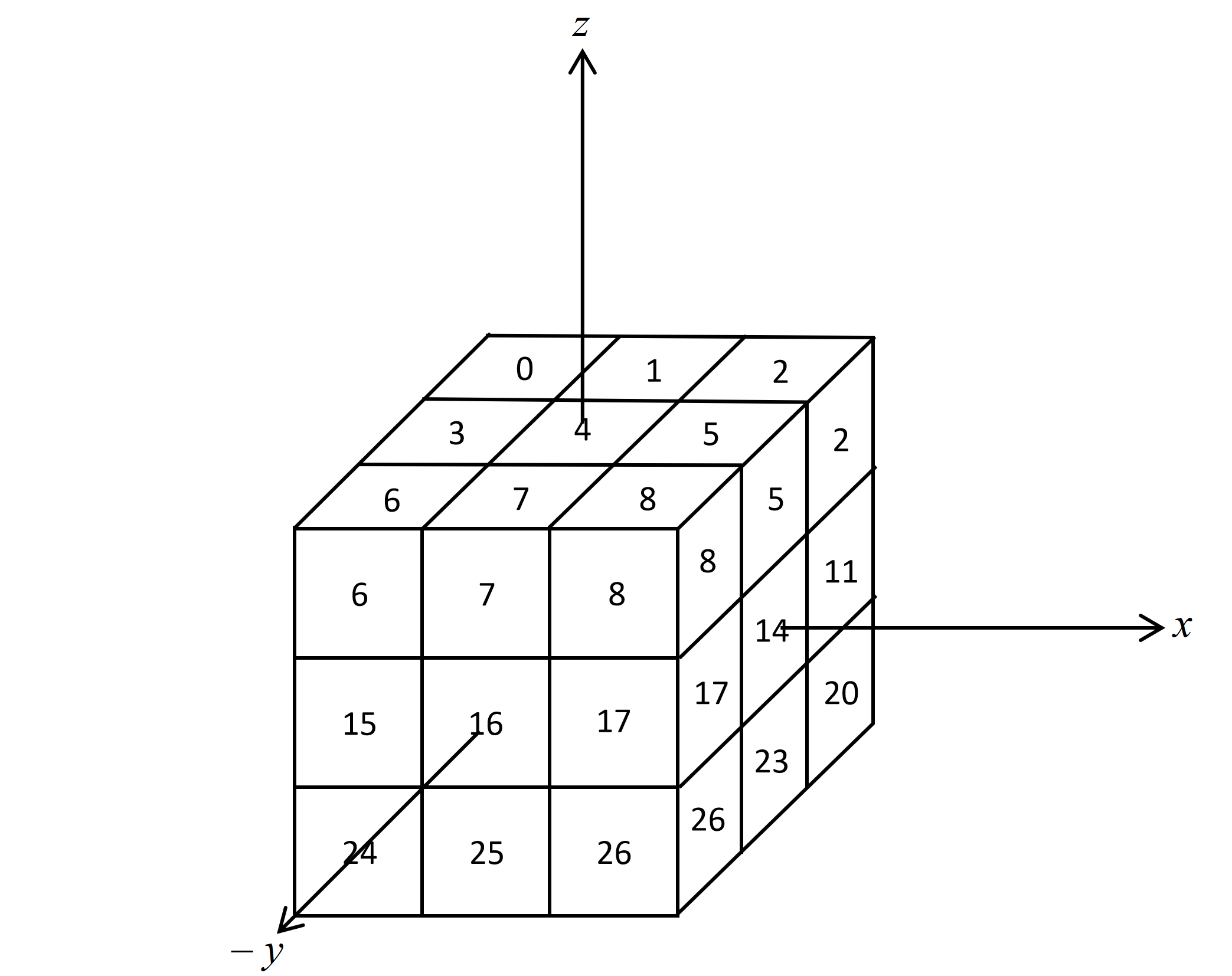}&
    \includegraphics[width=0.48\textwidth]{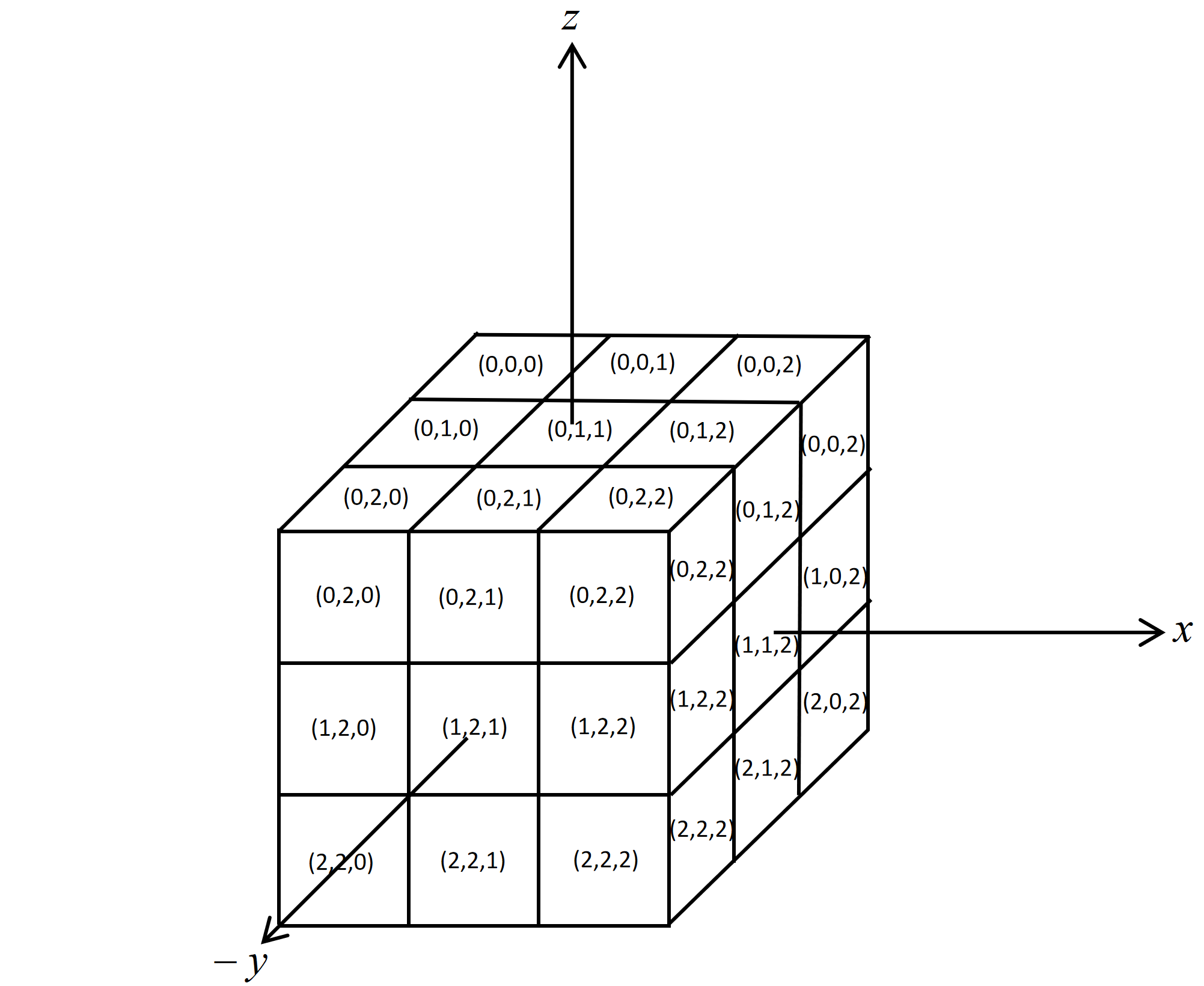}\\
    (a) & (b)
    \end{tabular}
    }
    \caption{(a) Voxel index in \ac{1D} form; (b) Voxel index in \ac{3D} form}
    \label{fig:Voxel_index_1D_Voxel_index_3D}
\end{figure}

\paragraph{\ac{3D} imaging coordinate system.} 
Let $(x, y, z)$ be the coordinate system. For a scanning \ac{3D} image, supposing that the positive parameters $L_x$, $L_y$, $L_z$ and $N_x$, $N_y$ $N_z$ are given, where $(L_x, L_y, L_z)$ and $(N_x, N_y, N_z)$ determine the side lengths of the domain and the size of the image along $x$-, $y$- and $z$-axis, respectively. Let $(d_x, d_y, d_z) := (L_x/N_x, L_y/N_y, L_z/N_z)$ be the side lengths or scales of the voxel along the corresponding axes. Without loss of generality, we assume that $L_x = L_y = L_z$, $N_x = N_y = N_z$, accordingly, $d_x = d_y = d_z$, and the center of the image domain is at the origin of the coordinate system as displayed in \cref{fig:Voxel_index_1D_Voxel_index_3D} (a). For simplicity, we further assume that the scales $d_x = d_y = d_z = 1$. If the scale is not unity, the real value of X-ray transform just equals to the scale multiplying that value for the case with unity scale. 

\paragraph{Voxel indexes.} 
We also define two different indexes for the voxels of the image. The one is given in \ac{1D} form as $I = 0, 1, \ldots, N_xN_yN_z - 1$, which is shown as the example with $N_x = N_y = N_z = 3$ in \cref{fig:Voxel_index_1D_Voxel_index_3D} (a). The other one is presented in \ac{3D} form as $(k, j, i)$ for $k = 0 , 1, \ldots, N_z - 1$, $j = 0, 1, \ldots, N_y - 1$, $i = 0, 1, \ldots, N_x - 1$, which is shown as the example also with $N_x = N_y = N_z = 3$ in \cref{fig:Voxel_index_1D_Voxel_index_3D} (b). The relation between them can be expressed as $I = kN_yN_x + jN_y + i$, which can be easily validated by the example in \cref{fig:Voxel_index_1D_Voxel_index_3D}. They also can convert between each other. 

Here the first and the second indexes are used to indicate the voxel positions in the associated vectorized image and the original image, respectively. Note that the second one is different from the coordinate position of the voxel. 

In this section, we will use voxel $(k, j, i)$ to indicate the voxel with index $(k, j, i)$. The valid $k$ (or $j$, $i$) is for itself satisfying $0 \leq k \leq N_z -1$ (or $0 \leq j \leq N_y -1$, $0 \leq i \leq N_x -1$). 

\paragraph{Eulerian angles.} 
To specify an arbitrary ray in \ac{3D} case, we need to use the concept of \textit{Eulerian angles} (see \cite[Chapter 5]{frank2006three}). Assuming that $(x, y ,z)$ is the original coordinate system. If the direction of a ray is $\theta$, defined by the three Eulerian angles $\phi_1$, $\phi_2$ and $\phi_3$, a projection is acquired on the plane that is perpendicular to the ray and also containing the origin. Then a new coordinate system is constructed by the above direction and plane, which is denoted by $(x', y', z')$ coordinate system. The transformation between the vectors in $(x, y ,z)$ coordinate system and those in $(x', y', z')$ coordinate system is given by three Eulerian rotations as 
\begin{equation}\label{eq:Eulerian_rotations}
\left[
   \begin{array}{c}
     x' \\
     y' \\
     z' \\ 
   \end{array}
\right] 
= 
\mathcal{R}
\left[
\begin{array}{c}
     x \\
     y \\
     z \\ 
   \end{array}
\right], 
\end{equation}
where 
\[
\mathcal{R} = \left[
   \begin{array}{ccc}
     \cos\phi_3 & \sin\phi_3 & 0 \\
     -\sin\phi_3 & \cos\phi_3 & 0 \\
     0 & 0 & 1 \\ 
   \end{array}
\right] \!\!
\left[
   \begin{array}{ccc}
     \cos\phi_2 & 0 & \sin\phi_2 \\
     0 & 1 & 0 \\ 
     -\sin\phi_2 & 0 & \cos\phi_2 \\
   \end{array}
\right] \!\!
\left[
   \begin{array}{ccc}
     \cos\phi_1 & \sin\phi_1 & 0 \\
     -\sin\phi_1 & \cos\phi_1 & 0 \\
     0 & 0 & 1 \\ 
   \end{array}
\right]. 
\]
These rotations can be illuminated as that first the $(x, y, z)$ coordinate system is contra-rotated by the angle $\phi_1$ around its $z$-axis, resulting in the intermediate coordinate system $(x_1, y_1,z_1)$, then by the angle $\phi_2$ around its new $y$-axis, yielding the second intermediate coordinate system $(x_2, y_2, z_2)$, and finally by the angle $\phi_3$ around its new $z$-axis to lead the final coordinate system $(x', y', z')$.

\subsection{Derivation of intersection length}
\label{sec:derivation_intersection_3D}

The purpose of this section is to derive intersection length for various \ac{3D} scanning geometries. 

\subsubsection{3D parallel beam}
\label{sec:3D_parallel_beam}

To characterize the \ac{3D} parallel beam, we can use the Eulerian angles above. Given that $(x, y ,z)$ is the coordinate system affixed to the image to be reconstructed. By projecting the image along a direction $\theta$, a \ac{2D} projection is acquired on its perpendicular plane containing the origin. Since the last rotation in \cref{eq:Eulerian_rotations} is an in-plane rotation, the parallel-beam scanning geometry can be generated by letting $\phi_3 = 0$. Then the transformation \cref{eq:Eulerian_rotations} becomes 
\begin{equation}\label{eq:Eulerian_rotations_2}
\left[
   \begin{array}{c}
     x' \\
     y' \\
     z' \\ 
   \end{array}
\right] 
= 
\left[
   \begin{array}{ccc}
     \cos\phi_2\cos\phi_1 & \cos\phi_2\sin\phi_1 & \sin\phi_2 \\
     - \sin\phi_1 & \cos\phi_1 & 0 \\ 
     - \sin\phi_2\cos\phi_1 & - \sin\phi_2\sin\phi_1 & \cos\phi_2 \\
   \end{array}
\right]\!\!
\left[
\begin{array}{c}
     x \\
     y \\
     z \\ 
   \end{array}
\right]. 
\end{equation}
Hence the direction of the parallel beam can be expressed as 
\begin{equation}\label{eq:3Dparallel_direction}
\theta = (\cos\phi_2\cos\phi_1, \cos\phi_2\sin\phi_1, \sin\phi_2), 
\end{equation}
and the associated orthogonal projection plane can be spanned by $\{\theta_1, \theta_2\}$, where 
\begin{align}
\label{eq:theta_1}&\theta_1 = (-\sin\phi_1, \cos\phi_1, 0), \\
\label{eq:theta_2}&\theta_2 = (-\sin\phi_2\cos\phi_1, -\sin\phi_2\sin\phi_1, \cos\phi_2). 
\end{align}

For the purpose of calculating the intersection length, the position of the ray and the pixel need to be specified. Assuming that the parameters $s_1, s_2 \in \Real^1$ and $\phi_1 \in [0, 2\pi)$, $\phi_2 \in [0, \pi)$ are given, where the $(s_1, s_2, \phi_1, \phi_2)$ determines a ray of the parallel beam in \ac{3D} space as shown in \cref{fig:3D_parallel_beam_circular_cone_beam_Helical_cone_beam} (a). Specifically, the $(s_1, s_2)$ is the coordinate position of the ray projecting onto the plane $\{\theta_1, \theta_2\}$, and the $(\phi_1, \phi_2)$ is the Eulerian angles indicating the direction $\theta$ of the ray. 

As assumed, the scales $d_x = d_y = d_z = 1$. We then give the transformation between the index of the voxel and the coordinate position of its center. For any voxel with index $(k, j, i)$, its center is located at the coordinate position 
\[
\biggl(i - \frac{N_x - 1}{2}, ~ \frac{N_y - 1}{2} - j,~\frac{N_z - 1}{2} - k\biggr).
\]
The associated voxel basis function is defined as 
\begin{equation}\label{eq:voxel_basis_function}
    f_{kji}(x, y, z) = \begin{cases}
  1, & (x, y, z) \in \Omega_{kji}, \\
  0, & \text{otherwise}.
 \end{cases}
\end{equation}
where $\Omega_{kji} := \bigl[i - \frac{N_x}{2}, i - \frac{N_x}{2} + 1\bigr] \times \bigl[\frac{N_y}{2} - j - 1, \frac{N_y}{2} - j\bigr] \times \bigl[\frac{N_z}{2} - k - 1, \frac{N_z}{2} - k\bigr]$ denotes the support of voxel $(k, j, i)$. 

The intersection length $l_{kji}(s_1, s_2, \phi_1, \phi_2)$ between the ray $(s_1, s_2, \phi_1, \phi_2)$ and the voxel $(k, j, i)$ can be expressed as the X-ray transform of the associated voxel basis function by 
\begin{equation}\label{eq:intersection_length_3D}
    l_{kji}(s_1, s_2, \phi_1, \phi_2) = \int_{-\infty}^{\infty}f_{kji}(t\theta + s_1\theta_1 + s_2\theta_2)\dint t.
\end{equation}
By \cref{eq:3Dparallel_direction}--\cref{eq:voxel_basis_function}, the \cref{eq:intersection_length_3D} implies that if and only if the following condition is satisfied  
\begin{equation}\label{eq:condition_nonzero_3D}
\begin{cases}
  i - \frac{N_x}{2} \leq t\cos\phi_2\cos\phi_1 - s_1\sin\phi_1 - s_2\sin\phi_2\cos\phi_1 \leq i - \frac{N_x}{2} + 1,  \\
  \frac{N_y}{2} - j - 1 \leq t\cos\phi_2\sin\phi_1 + s_1\cos\phi_1 - s_2\sin\phi_2\sin\phi_1 \leq \frac{N_y}{2} - j, \\
   \frac{N_z}{2} - k - 1 \leq t\sin\phi_2 + s_2\cos\phi_2 \leq \frac{N_z}{2} - k, 
 \end{cases}
\end{equation}
then the length of intersection might be non-vanishing. Note that the variable $t$ of integration in \cref{eq:intersection_length_3D} is along the ray. Therefore, the range of $t$ that satisfies the inequalities in \cref{eq:condition_nonzero_3D} would be the length of intersection. 

As we observe, the ray $(s_1, s_2, \phi_1 + 2n_1\pi, \phi_2 + 2n_2\pi)$ for $n_1, n_2 \in \mathbb{Z}$ is equivalent to the ray $(s_1, s_2, \phi_1, \phi_2)$, the ray $(s_1, s_2, \phi_1, \phi_2 + \pi)$ is the same as the ray $(s_1, -s_2, \phi_1, \phi_2)$, and the ray $(-s_1, s_2, \phi_1 + \pi, -\phi_2)$ is the same as the ray $(s_1, s_2, \phi_1, \phi_2)$. Conclusively, any ray can be specified by the ray $(s_1, s_2, \phi_1, \phi_2) \in \Real^1 \times \Real^1 \times [0, 2\pi) \times [0, \pi/2]$. By the analysis above, the possible cases has been largely reduced. In what follows we list the cases by the different choices of $(\phi_1, \phi_2)$. 

\paragraph{Case 1: $(\phi_1, \phi_2) \in (0, \pi/2) \times (0, \pi/2)$.} For simplicity, let 
\begin{align}\label{eq:C_xC_yC_z}
\nonumber &C_x = \frac{i - N_x/2 + s_1\sin\phi_1 + s_2\sin\phi_2\cos\phi_1}{\cos\phi_2\cos\phi_1}, \\
&C_y = \frac{N_y/2 - j - s_1\cos\phi_1 + s_2\sin\phi_2\sin\phi_1}{\cos\phi_2\sin\phi_1}, \quad 
C_z = \frac{N_z/2 - k - s_2\cos\phi_2}{\sin\phi_2}. 
\end{align}
By simple calculations, using \cref{eq:C_xC_yC_z}, the condition \cref{eq:condition_nonzero_3D} becomes  
\begin{equation}\label{eq:condition_nonzero_t_3D}
\begin{cases}
  C_x \leq t  \leq C_x + \frac{1}{\cos\phi_2\cos\phi_1},  \\
  C_y - \frac{1}{\cos\phi_2\sin\phi_1} \leq t \leq C_y,  \\
  C_z - \frac{1}{\sin\phi_2} \leq t \leq C_z. 
 \end{cases}
\end{equation}

Obviously, if the three intervals in \cref{eq:condition_nonzero_t_3D} have no overlap or their overlap is just a point, then the length of the intersection is absolutely zero. To rule out this possibility, if and only if 
\begin{equation}\label{eq:conditions_3D_1}
C_{low} < C_{up}. 
\end{equation}
Here 
\begin{align*}
C_{low} &= \max\bigl(C_x, C_y - 1/(\cos\phi_2\sin\phi_1), C_z - 1/\sin\phi_2\bigr), \\
C_{up} &= \min(C_x + 1/(\cos\phi_2\cos\phi_1), C_y, C_z).  
\end{align*}
Equivalently, \cref{eq:conditions_3D_1} can be rewritten as 
\begin{equation}\label{eq:conditions_3D_2}
\begin{cases}
 C_x < C_y, \quad C_y - \frac{1}{\cos\phi_2\sin\phi_1} < C_x + \frac{1}{\cos\phi_2\cos\phi_1}, \\
 C_x < C_z, \quad C_z - \frac{1}{\sin\phi_2} < C_x + \frac{1}{\cos\phi_2\cos\phi_1}, \\
 C_y - \frac{1}{\cos\phi_2\sin\phi_1} < C_z,  \quad C_z - \frac{1}{\sin\phi_2} < C_y. 
 \end{cases}
\end{equation}
For this case, the \cref{eq:conditions_3D_1} or \cref{eq:conditions_3D_2} is exactly the \textit{sufficient and necessary condition} for non-vanishing intersectability of the given ray and voxel. 

More precisely, by \cref{eq:C_xC_yC_z}, if the $i$ is given, the first two inequalities of \cref{eq:conditions_3D_2} reads as the following inequality 
\begin{equation}\label{eq:restriction_j_3D}
C_1  - \tan\phi_1 - 1 < j < C_1. 
\end{equation}
Here 
\[
        C_1 = \frac{N_y}{2} - s_1\cos\phi_1 + s_2\sin\phi_2\sin\phi_1 - C_x\cos\phi_2\sin\phi_1. 
\]
Fixed $i$, once we obtain the range of valid $j$ by \cref{eq:restriction_j_3D}, the range of valid $k$ can be calculated by the middle and last two inequalities of \cref{eq:conditions_3D_2} as 
\begin{equation}\label{eq:restriction_k_3D_1}
\begin{cases}
C_2  - \frac{\tan\phi_2}{\cos\phi_1} - 1 < k < C_2, \\
C_3 - 1 < k < C_3 + \frac{\tan\phi_2}{\sin\phi_1}. 
\end{cases}
\end{equation}
Here  
\[
        C_2 = \frac{N_z}{2} - s_2\cos\phi_2 - C_x\sin\phi_2, \quad C_3 = \frac{N_z}{2} - s_2\cos\phi_2 - C_y\sin\phi_2.
\]

Hence for any given $i$, we can calculate the range of valid $j$ by \cref{eq:restriction_j_3D}, and further that of $k$ by the two inequalities in \cref{eq:restriction_k_3D_1}, to obtain those voxels intersecting with the given ray non-vanishingly. Similarly, provided that the $j$ or $k$ is given, the valid range of $i$, $k$ or $i$, $j$ can be calculated by the same method above. For those voxels, the condition \cref{eq:conditions_3D_1} is definitely fulfilled. By merging the intervals in \cref{eq:condition_nonzero_t_3D}, the intersection reads as 
\begin{equation}\label{eq:range_t_3D}
   C_{low} \leq t \leq C_{up}. 
\end{equation} 
Otherwise, the intersection is empty. Hence, the \textit{analytic formula} of the non-vanishing intersection length would be 
\begin{equation}\label{eq:range_length_3D}
    l_{kji}(s_1, s_2, \phi_1, \phi_2) = C_{up} - C_{low}.
\end{equation}

Specifically, we give a sketch of the algorithm for $(\phi_1, \phi_2) \in (0, \pi/2) \times (0, \pi/2)$ as follows. 
\begin{algorithm}
\caption{A sketch of the algorithm for $(\phi_1, \phi_2) \in (0, \pi/2) \times (0, \pi/2)$. }
\label{algo:algorithm_3D}
\begin{algorithmic}[1]
\State Given $s_1, s_2, \phi_1, \phi_2$ to determine a ray.
\State \emph{Loop}: For $0 \leq i \leq N_x - 1$ 
\State \quad\quad\quad\quad\quad Obtain the valid $j$ by \cref{eq:restriction_j_3D}. 
\State \quad\quad\quad\quad\quad Compute the valid $k$ by \cref{eq:restriction_k_3D_1}. 
\State \quad\quad\quad\quad\quad Calculate the intersection length by \cref{eq:range_length_3D} for voxel ($k, j, i$). 
\State \emph{End loop}.
\State \textbf{Output} the intersection lengths and associated \ac{1D}-form indexes.
\end{algorithmic}
\end{algorithm}

\begin{remark}\label{rem:case1_3D}
For other similar cases, such as $(\phi_1, \phi_2) \in \{(\pi/2, \pi)\cup(\pi, 3\pi/2)\cup(3\pi/2, 2\pi)\} \times (0, \pi/2)$, the corresponding sufficient and necessary condition for non-vanishing intersectability, and analytic formula of the intersection length can be easily obtained by making quite minor changes for the derivation above. So we skip these negligible derivations. 
\end{remark}

\paragraph{Case 2: $(\phi_1, \phi_2) = \{0\} \times (0, \pi/2)$.}  
Using simple computations, the condition \cref{eq:condition_nonzero_3D} is translated into 
 \begin{equation}\label{eq:condition_nonzero_t_case2_3D_1}
\begin{cases}
 C_x \leq t \leq C_x + \frac{1}{\cos\phi_2}, \quad C_z - \frac{1}{\sin\phi_2} \leq t \leq C_z, \\
  \frac{N_y}{2} - j - 1 \leq s_1 \leq \frac{N_y}{2} - j.
 \end{cases}
\end{equation}
Here 
\begin{equation}\label{eq:C_xC_z}
C_x = \frac{ i - N_x/2 + s_2\sin\phi_2}{\cos\phi_2},  \quad C_z = \frac{N_z/2 - k - s_2\cos\phi_2}{\sin\phi_2}.
\end{equation}
Following the derivation in \textbf{case 1} of this section, if all the valid $j$ are outside the interval $\frac{N_y}{2} - s_1 - 1 \leq j \leq \frac{N_y}{2} - s_1$, or the first and last intervals in \cref{eq:condition_nonzero_t_case2_3D_1} have no overlap or their overlap is only a point, the length of the intersection is absolutely vanishing. To exclude this possibility, if and only if 
\begin{equation}\label{eq:conditions_3D_case2_1}
\begin{cases}
C_{low} < C_{up}, \\
\frac{N_y}{2} - s_1 - 1 \leq j \leq \frac{N_y}{2} - s_1.
\end{cases} 
\end{equation}
Here 
\[
C_{low} = \max\bigl(C_x, C_z - 1/\sin\phi_2\bigr), \quad C_{up} = \min(C_x + 1/\cos\phi_2, C_z).  
\]
The \cref{eq:conditions_3D_case2_1} is equivalent to 
\begin{equation}\label{eq:condition_nonzero_t_case2_3D}
\begin{cases}
 C_x < C_z, \quad C_z - \frac{1}{\sin\phi_2} < C_x + \frac{1}{\cos\phi_2}, \\
 \frac{N_y}{2} - s_1 - 1 \leq j \leq \frac{N_y}{2} - s_1. 
 \end{cases}
\end{equation}
For this case, the \cref{eq:conditions_3D_case2_1} or \cref{eq:condition_nonzero_t_case2_3D} is the \textit{sufficient and necessary condition} for non-vanishing intersectability of the given ray and voxel. 

Moreover, by \cref{eq:C_xC_z}, if the $i$ is given, the first two inequalities of \cref{eq:condition_nonzero_t_case2_3D} becomes the following inequality 
\begin{equation}\label{eq:restriction_k_case2_3D}
C_1  - \tan\phi_2 - 1 < k < C_1. 
\end{equation}
Here 
\[
        C_1 = \frac{N_z}{2} - s_2\cos\phi_2 - C_x\sin\phi_2. 
\]
If the valid $j$ satisfies $\frac{N_y}{2} - s_1 - 1 \leq j \leq \frac{N_y}{2} - s_1$, for given $i$, we can calculate the range of valid $k$ by \cref{eq:restriction_k_case2_3D},   
to obtain those voxels intersecting with the given ray non-vanishingly. The \textit{analytic formula} of the non-vanishing intersection length would be 
\begin{equation}\label{eq:range_length_case2_3D}
    l_{kji}(s_1, s_2, \phi_1, \phi_2) = C_{up} - C_{low}.
\end{equation}

\begin{remark}\label{rem:case2_3D}
For such cases $(\phi_1, \phi_2) = \{\pi/2, \pi, 3\pi/2\} \times (0, \pi/2)$, the sufficient and necessary condition for non-vanishing intersectability, and analytic formula of the intersection length can be immediately computed following the method in this case. So the trivial derivations are omitted here. Note that the rays in such cases are actually parallel to a certain coordinate plane, which leads the corresponding problem reduce to a \ac{2D} problem as derived above. However, the ambiguity perhaps happen (see the third example in \cref{sec:ambiguity_solution}). 
\end{remark}

\paragraph{Case 3: $(\phi_1, \phi_2) = (0, 0)$.} 
With simple calculations, the condition \cref{eq:condition_nonzero_3D} reads as  
 \begin{equation}\label{eq:condition_nonzero_t_case3_3D}
\begin{cases}
  i - \frac{N_x}{2} \leq t \leq i - \frac{N_x}{2} + 1,  \\
  \frac{N_y}{2} - j - 1 \leq s_1 \leq \frac{N_y}{2} - j, \\
  \frac{N_z}{2} - k - 1 \leq s_2 \leq \frac{N_z}{2} - k.  
 \end{cases}
\end{equation}
Evidently, for any given $i$, if the valid $j$ satisfies $\frac{N_y}{2} - s_1 - 1 \leq j \leq \frac{N_y}{2} - s_1$, and the valid $k$ also satisfies $\frac{N_z}{2} - s_2 - 1 \leq k \leq \frac{N_z}{2} - s_2$, the intersection length is unity, otherwise, it is zero.

\begin{remark}\label{rem:case3_3D}
For the cases as $(\phi_1, \phi_2) \in \{\{\pi/2, \pi, 3\pi/2\} \times \{0\}\}$, the intersection length can be computed by the method of the case above. So we omit the details. 
\end{remark}

\paragraph{Case 4: $(\phi_1, \phi_2) = [0, 2\pi) \times \{\pi/2\}$.} 
With simple calculations, the condition \cref{eq:condition_nonzero_3D} is rewritten as 
 \begin{equation}\label{eq:condition_nonzero_t_case4_3D}
  \begin{cases}
  i - \frac{N_x}{2} \leq - s_1\sin\phi_1 - s_2\cos\phi_1 \leq i - \frac{N_x}{2} + 1,  \\
  \frac{N_y}{2} - j - 1 \leq s_1\cos\phi_1 - s_2\sin\phi_1 \leq \frac{N_y}{2} - j, \\
   \frac{N_z}{2} - k - 1 \leq t  \leq \frac{N_z}{2} - k. 
 \end{cases}
\end{equation}
Hence, for any given $k$, if the valid $i$ satisfies $\frac{N_x}{2} - s_1\sin\phi_1 - s_2\cos\phi_1 - 1 \leq i \leq \frac{N_x}{2} - s_1\sin\phi_1 - s_2\cos\phi_1$, and the valid $j$ also satisfies $\frac{N_y}{2} - s_1\cos\phi_1 + s_2\sin\phi_1 - 1 \leq j \leq \frac{N_y}{2} - s_1\cos\phi_1 + s_2\sin\phi_1$, the intersection length is unity, otherwise, it is zero. 

\begin{remark}\label{rem:case4_3D}
Obviously, the ray in \textbf{cases 3} and \textbf{4} is parallel to a certain axis, which results in quite simple calculations. But the ambiguity perhaps happen (see the second example in \cref{sec:ambiguity_solution}). 
\end{remark}

\subsubsection{Circular cone beam}
\label{sec:3D_cone_beam}

As mentioned previously, the cone beam is a projection geometry that is often used in \ac{3D} circumstance. Next, we will generalize the algorithm to this situation. 

\paragraph{Equiangular circular cone beam.} 
The equiangular cone beam is an extension of the equiangular fan beam into \ac{3D} situation. Actually, any ray of the equiangular cone beam can be specified by the parameters $(D, \phi'_1, \alpha, \beta)$, as illustrated in \cref{fig:3D_parallel_beam_circular_cone_beam_Helical_cone_beam} (b). The source locates at $xOy$ plane. Here $D > 0$ indicates the distance of source for X-rays from the origin, $\phi'_1 \in [0, 2\pi)$ stands for the angle between the center line from the source to origin and the positive $x$-axis, and $\alpha \in (-\pi/2, \pi/2)$ determines the signed angle between the projection of the ray on $xOy$ plane and the center line, and $\beta \in (-\pi/2, \pi/2)$ denotes the signed angle between the ray and its projection on $xOy$ plane. 

Once given the parameters $(D, \phi'_1, \alpha, \beta)$, we can convert them into the parameters used in \ac{3D} parallel beam by performing the simple transformation as the following  
\begin{equation}\label{eq:equiangular_cone_beam}
        \phi_1 = \phi'_1 + \alpha, \quad
        \phi_2 = \beta, \quad
        s_1 = D\sin\alpha, \quad
        s_2 = D\cos\alpha\sin\beta. 
\end{equation}

\paragraph{Equispaced circular cone beam.} 
The equispaced circular cone beam can also be derived in a similar way. The only difference is that it locates the ray using two distances $(t, h)$ rather than two angles $(\alpha, \beta)$, 
as shown in \cref{fig:3D_parallel_beam_circular_cone_beam_Helical_cone_beam} (b). The two distances specify the coordinate position of the detector in the projection plane. In other words, the parameters for equispaced cone beam is $(D, \phi'_1, t, h)$. Similarly, we can convert them into the 
parameters in \ac{3D} parallel beam as 
\begin{equation}\label{eq:equispaced_cone_beam}
        \phi_1 = \phi'_1 + \alpha, \quad
        \phi_2 = \beta, \quad
        s_1 = D\sin\alpha, \quad
        s_2 = h\cos^2\alpha\cos\beta. 
\end{equation}
Here  
\[
\alpha = \arctan\frac{t}{D}, \quad \beta = \arctan\biggl(\frac{h}{\sqrt{D^2 + t^2}}\biggr).
\]

After transformation, the associated X-ray transform in circular cone beam can be readily computed by the algorithm developed in \cref{sec:3D_parallel_beam}.

\subsubsection{Helical cone beam}
\label{sec:Helical_cone_beam}

Helical cone beam is actually a type of cone beam that swipes through the detected body while the travelling path of its X-ray source relative to a fixed point on the body is shaped like a helix, as depicted in \cref{fig:3D_parallel_beam_circular_cone_beam_Helical_cone_beam} (c). Hence the source would have an additional parameter, i.e., the signed vertical distance to the origin. 

\paragraph{Equiangular helical cone beam.} 
Here the parameters are $(D, \phi'_1, \beta_1, \beta_2, H)$, where $H$ determines the signed vertical distance from source to origin. Since in \cref{sec:3D_cone_beam} only $s_2$ is related to this vertical distance, to obtain the transformation, we just need to modify the last formula in \cref{eq:equiangular_cone_beam} as  
\begin{equation}\label{eq:Helical_cone_beam_equiangular}
                s_2 = D\cos\alpha\sin\beta + H\cos\beta.
\end{equation} 
The other formulas are still unchanged. 

\paragraph{Equispaced helical cone beam.} 
Similarly, the parameters for this beam is $(D, \phi'_1, t, h, H)$. To obtain the conversion, we only change the last formula of \cref{eq:equispaced_cone_beam} as  
\begin{equation}\label{eq:Helical_cone_beam_equispaced}
       s_2 = h\cos^2\alpha\cos\beta + H\cos\beta, 
\end{equation}
and maintain the other formulas. 

After minor changes, the X-ray transform in helical cone beam can be computed by the algorithm developed in \cref{sec:3D_parallel_beam}. 

\begin{remark}\label{rem:adaptive_geometry}
The proposed algorithm is not just limited to the above circular/helical cone beams, which is applicable to any beam that its parameters can be transformed into those of the parallel beam. 
\end{remark}

\section{Discussion}
\label{sec:Discussions}

For the problem itself, the ambiguities would happen to some special situations. This is intrinsic to the computation of the X-ray transform. Here we will have a discussion, and give a solution in the proposed algorithm. And the adaptability and parallelization of the algorithm will be further discussed. 

\subsection{Ambiguities and solution}
\label{sec:ambiguity_solution}

\paragraph{Ambiguities.} 
For the kind of situations as \textbf{case 2} in \cref{sec:2D_parallel_beam} and \textbf{cases 3} and \textbf{4} in \cref{sec:3D_parallel_beam}, the given ray is parallel to a certain axis, and the calculation of intersection is quite simple. More precisely, the length of intersection is just unity if it is non-vanishing. Nevertheless, when the given ray is exactly overlapping with grid line in \ac{2D}/\ac{3D} scenarios, the ambiguity would happen to the problem itself. 

As an example in \ac{2D} situation, let $N_x = N_y = 5$, and take the ray $(1.5, 0)$ (i.e., $s = 1.5$, $\phi = 0$). By simple computations, for any given $i$, the intersection length is unity if $0 \leq j \leq 1$. Obviously, the same one intersection can be attributed to the ray $(1.5, 0)$ with two different pixels $(0, i)$ and $(1, i)$ for $0 \leq i \leq 4$, simultaneously, which results in ambiguity. 

Another example is given in \ac{3D} case. Let $N_x = N_y = N_z = 4$, and consider the ray $(1, 1, 0, 0)$ (i.e., $s_1 = 1$, $s_2 = 1$, $\phi_1 = 0$, $\phi_2 = 0$). For each $i$, if $0 \leq j \leq 1$ and $0 \leq k \leq 1$, the intersection length is unity. That is to say, the only one intersection can be assigned to the ray $(1, 1, 0, 0)$ with four different voxels $(0, 0, i)$, $(0, 1, i)$, $(1, 0, i)$ and $(1, 1, i)$ for $0 \leq i \leq 3$ coincidently, which also leads to the ambiguity. 

The other ambiguity would occur when the given ray is accurately on the grid plane, which is a special situation as described in \textbf{case 2} for \ac{3D} circumstance. For instance, let $N_x = N_y = N_z = 3$, think about the ray $(-0.5, \sqrt{2}, 0, \pi/4)$ (i.e., $s_1 = -0.5$, $s_2 = \sqrt{2}$, $\phi_1 = 0$, $\phi_2 = \pi/4$). By simple calculations, we found that $1 \leq j \leq 2$, and only for $i = 0$, there exists the valid $k = 0$, and the length of intersection is $\sqrt{2}$, which is consistent with the geometric observation. In other words, the only one intersection can be assigned to the ray $(-0.5, \sqrt{2}, 0, \pi/4)$ with two different voxels $(0, 1, 0)$ and $(0, 2, 0)$ simultaneously, which produces the ambiguity. 

\paragraph{Solution.} 
If the ambiguities above appear, to get rid of them, we can always consistently take the pixel/voxel with bigger (or lower) \ac{1D}-form index (see \cref{fig:2Dindex_2Dindex2Dform} (a) and \cref{fig:Voxel_index_1D_Voxel_index_3D} (a)) as the intersected unit with the given ray in the proposed algorithm. For the first example, that is the pixel $(1, i)$. Because the pixel/voxel index is independent of the coordinate system. In this way, for any ray overlapping with a grid line or being on a grid plane, the pixels/voxels on the same side of the ray would be chosen reasonably.

\subsection{Adaptability}
\label{sec:Adaptability}

The proposed algorithm is entirely obtained by mathematical derivations, which needs not compute any intersection point of the ray and the grid lines or planes. It is worth noting that even if the center of the image is not at the origin of the imaging coordinate system, or the scales $(d_x, d_y, d_z)$ of the pixel/voxel or the sizes $(N_x, N_y, N_z)$ of the image are different from each other along associated axes, the algorithm can be derived by the same pipeline. 

The algorithm is independent of the specific scanning geometry, and can output the projection matrix in sparse storage. Because its fundamental element is to compute the intersection length of a certain ray with the pixels/voxels non-vanishingly intersecting with it. The algorithm can be immediately extended to deal with the scanning geometry even that the rays are random distributed. No matter how are the detectors distributed, the geometric parameters of the ray (see \cref{sec:2D_parallel_beam} and \cref{sec:3D_parallel_beam}) can be determined. Then the proposed algorithm can be used to compute the associated X-ray transform, which depends only on the geometric parameters of the ray. Hence, the algorithm can be customized freely according to the requirements of the users, and more scanning geometries can be easily added into the framework based on the proposed algorithm. 

Moreover, the proposed algorithm is also suited to the case with polygonal/polyhedral basis function, whose special forms include pixel/voxel basis function. As we know, the support of polygonal/polyhedral basis function is a simplex. Along the method in \cref{sec:2D_parallel_beam} and \cref{sec:3D_parallel_beam}, it is easy to obtain the sufficient and necessary condition for non-vanishing intersectability, and the analytic formula for the intersection length.

\subsection{Computational cost and parallelization}
\label{sec:Computational_cost}

By the condition of non-vanishing intersectability, we can distinguish out the units that produce valid intersections with the given ray. Only for these units rather than all the individuals, we then calculate the intersection lengths by the obtained analytic formula. More precisely, for any given $i$ (or $j$), we can calculate the range of valid $j$ (or $i$) by the condition \cref{eq:restriction_j} in \ac{2D} situation, to obtain those pixels intersecting with the given ray non-vanishingly. In \ac{3D} situation, for any given $i$ (or $j$, $k$), we can calculate the ranges of valid $j$ and $k$ (or $i$ and $k$, $i$ and $j$) by the conditions \cref{eq:restriction_j_3D} and \cref{eq:restriction_k_3D_1}, to obtain those voxels intersecting with the given ray non-vanishingly. And then we just calculate the intersection lengths by the obtained analytic formula for those units rather than all the individuals. The overhead that comes from evaluating the ranges is at most one addition and one multiplication for \ac{2D} case, and triple additions and triple multiplications for \ac{3D} case. The computational complexity is $O(1)$ from computing the analytic formula in \cref{eq:range_length} or \cref{eq:range_length_3D}. So the computational cost of the proposed algorithm is $O(N)$ for any given ray, and $O(NM)$ for all rays even to \ac{3D} circumstance. Here $N$ represents the size of the image being reconstructed along one axis, and $M$ denotes the number of rays. Hence, this algorithm achieves the degree of optimality, and its computational cost is optimal. 
 
Since the proposed algorithm implements the calculation of intersection length ray by ray and unit by unit, and all of the calculations are independent of each other, it is quite suited to parallelize and the computational complexity per parallel thread achieves $O(1)$.

\subsection{Comparison with existing approaches}
\label{sec:Comparison}

In this part, the brief comparison will be made among the proposed model and several mainly existing approaches (i.e., Siddon's method in \cite{siddon_1985} and its improved method in \cite{gao_2012}). 

As mentioned previously, Siddon's method requires to compute the intersection points of the ray with all of the grid lines/planes in 2D/3D circumstances, and then sorts all those intersection points. This approach is time-consuming, its computational cost is proportional to the size of the image, and it is impossible to parallelize the sorting procedure. The improved version in \cite{gao_2012} always needs to compute the intersection of a ray with two end points with a 2D/3D domain of the image, and then uses the relation of the two end points to classify the cases of intersection of the ray with coordinate axes/planes. Actually, it is mildly time-consuming to compute the two end points for each ray, where the method to compute the end points is not covered in \cite{gao_2012}. 

However, the proposed method is quite different from the above existing methods, which eliminates the computational process of the end points, but classifies the cases by the parameters of the ray directly. Hence, this algorithm is more intuitive and faster than the above approaches for computing the X-ray transform. And the ambiguity and its solution are also discussed in \cref{sec:ambiguity_solution}, which is necessary for the problem itself. However, these issues are not involved in \cite{siddon_1985} and \cite{gao_2012}. Hence, the proposed algorithm is more complete to use. 

Moreover, the proposed algorithm can be applied to deal with various 2D/3D scanning geometries straightforwardly, particularly for the complicatedly but commonly used 3D circular and helical cone beams, and even more general scanning geometries, since we have derived the associated transformation formulas explicitly. But these topics are not discussed in \cite{siddon_1985} and \cite{gao_2012}. As discussed in \cref{sec:Adaptability}, the proposed algorithm is suited to the general case with polygonal/polyhedral basis function. But the algorithms in \cite{siddon_1985} and \cite{gao_2012} depend strictly on the special structure of pixel/voxel, which are not applicable to the general case above. Hence, the proposed algorithm is more flexible, e.g., one can easily work with different scanning geometries and different image basis functions.

\section{Validations}
\label{sec:numerical_tests}

The algorithm proposed above has been validated through different test examples. Specific tests can be found as follows. This algorithm calculates the X-ray transform accurately using mathematical derivations, so the numerical errors do not exist other than the machine error. However, since the calculation is done in float data type and the program is written in C++ language, the result only has six significant digits. Since the scales of each unit (pixel/voxel) are fixed to be unity, the output is accurate down to around five decimal places, which should be accurate enough in most cases. If the higher precision is required, the program can calculate down to $15$ decimal places by switching float data type into double data type. Due to the limitation of space, not all tests can be given here. Thus, we will give only one test for each scanning geometry. Although this is not a full evaluation, the tests illustrate the correctness of the proposed method.

\subsection{Test suite 1: \ac{2D} parallel beam} 
\label{sec:Test_suite_1}

In \ac{2D} circumstance, the first example is for \ac{2D} parallel beam. Assuming that the size of the image is $3\times 3$. In this situation, one easy-to-understand example is a ray with $\pi/4$ angle to the positive $x$-axis and a distance of unity to the origin, namely, the ray $(1, \pi/4)$. 

Using intuitively geometric computation, the valid intersections of the ray with all the pixels can be fast obtained, where the intersected pixels are $(0, 0)$, $(0, 1)$ and $(1, 0)$, and the corresponding intersection lengths are $2 - \sqrt{2}$, $2\sqrt{2} - 2$ and $2\sqrt{2} - 2$, respectively.
    
On the other hand, we compute the results by hand along the algorithm in \cref{sec:2D_parallel_beam}. \\
(i) If $i = 0$, then $C_x = -3\sqrt{2}/2 + 1$, $C_1 = 3 - \sqrt{2}$, and then the valid $j$ includes $0$ and $1$. The $C_y = 3\sqrt{2}/2 - 1$ and $\sqrt{2}/2 - 1$ for $j = 0$ and $1$ respectively. Hence, the lengths of intersections of the ray with pixels $(0, 0)$ and $(1, 0)$ are $2 - \sqrt{2}$ and $2\sqrt{2} - 2$, respectively; \\
(ii) If $i = 1$, then $C_x = -\sqrt{2}/2 + 1$, $C_1 = 5/2 - 3\sqrt{2}/2$, and the valid $j$ only contains $0$. So the $C_y = 3\sqrt{2}/2 - 1$.  The length of intersection of the ray with pixel $(0, 1)$ is $2\sqrt{2} - 2$; \\
(iii) If $i = 2$, the valid $j$ is empty. 

By numerical implementation, the output of the valid lengths of intersections is shown in \cref{tab:test_2D_parallel}, with the first column being the \ac{1D}-form index of the pixel and the corresponding second column showing the intersection length. Keep in mind that the \ac{1D}-form index and the \ac{2D}-form one can be converted into each other.  
  \begin{table}[htbp]
  \centering
   \begin{tabular}{cl}
   \hline
    0 & 0.585787 \\ [-0.3em]
    1 & 0.828427 \\ [-0.3em]
    3 & 0.828427 \\ 
    \hline
    \end{tabular}
    \caption{Output of the example for \ac{2D} parallel beam}
    \label{tab:test_2D_parallel}
    \end{table}

Obviously, the above results are completely consistent with each other if we do not consider the machine error. In this way, the error is zero.

\subsection{Test suite 2: Fan beam} 
\label{sec:Test_suite_2}

Here we only give an example for the equiangular fan beam because the validation has no essential distinction with that of the equispaced one. For equiangular fan beam, the main algorithm are the same as the \ac{2D} parallel beam. What we just need to do is to get the corresponding parameters $(s, \phi)$ in \ac{2D} parallel beam for each ray when given $(D, \alpha, \gamma)$, as described in \cref{sec:fan_beam}. 

As an example, the size of the image is assumed to be $4\times 4$, the ray is considered with parameters $D = 4$, $\alpha = \pi/2$, $\gamma = -\pi/6$. With \cref{eq:2Dparalleltofan}, the parameters are transformed into the associated ones in \ac{2D} parallel beam as $s = -2$, $\phi = -\pi/6$. As observed in \cref{sec:2D_parallel_beam}, it indicates the ray $(2, 5\pi/6)$. 

By directly analytic geometry, the valid intersections of the ray with all the pixels can be gained immediately, where the intersected pixels are $(3, 0)$ and $(3, 1)$, and the corresponding lengths of intersections are $2\sqrt{3}/3$ and $4 - 2\sqrt{3}$, respectively. 

Then we recompute the results by hand along the algorithm in \cref{sec:2D_parallel_beam}. \\
(i) If $i = 0$, then $C_x = 2\sqrt{3}/3$, and then the valid $j$ should satisfy $1 + 2\sqrt{3}/3 < j < 2 + \sqrt{3}$. Hence, the valid $j$ only contains $3$. The $C_y = 2\sqrt{3} - 2$ for $j = 3$. Hence, the length of intersection of the ray with pixels $(3, 0)$ is $2\sqrt{3}/3$; \\
(ii) If $i = 1$, then $C_x = 0$, and the valid $j$ should satisfy $1 + \sqrt{3} < j < 2 + 4\sqrt{3}/3$. So the valid $j$ only contains $3$. Then $C_y = 2\sqrt{3} - 2$ as well. The length of intersection of the ray with pixel $(3, 1)$ is $4 - 2\sqrt{3}$; \\
(iii) If $i = 2$ or $3$, the valid $j$ is empty. 

By numerical implementation, the output of the valid lengths of intersections is provided in \cref{tab:test_fan}. The $12$th and $13$th pixels in \ac{1D} form are the pixels $(3, 0)$ and $(3, 1)$, respectively.  
  \begin{table}[htbp]
  \centering
   \begin{tabular}{cl}
   \hline
    12 & 1.1547 \\ [-0.3em]
    13 & 0.535899 \\ 
     \hline
    \end{tabular}
    \caption{Output of the example for equiangular fan beam}
    \label{tab:test_fan}
    \end{table}

It is easy to observe that the above results are completely consistent with each other if we do not take care the machine error. So the error is zero.

\subsection{Test suite 3: 3D parallel beam} 
\label{sec:Test_suite_3}

In \ac{3D} situation, we first consider the \ac{3D} parallel beam. Assuming that the size of the image is $3\times 3$. In this case, one intuitive example is a ray with $\pi/4$ angle to both the positive $x$- and $y$-axis and a distance of zero to the origin, namely, the ray $(0, 0, \pi/4, \pi/4)$. 

Not that the given ray locates at one diagonal plane of the image. We can immediately calculate out the valid intersections of the ray with all the voxels by geometric observation. The valid intersected voxels are $(0, 0, 2)$, $(0, 1, 1)$, $(1, 1, 1)$, $(2, 1, 1)$ and $(2, 2, 0)$, and the corresponding lengths of intersections are $3\sqrt{2}/2 -1$, $1 - \sqrt{2}/2$, $\sqrt{2}$, $1 - \sqrt{2}/2$ and $3\sqrt{2}/2 -1$, respectively. 

Then we recompute the results by hand along the algorithm in \cref{sec:3D_parallel_beam}. \\
(i) If $i = 0$, then $C_x = -3$, and then the valid $j$ should satisfy $1 < j < 3$. Hence, the valid $j$ only contains $2$. The $C_y = -1$ for $j = 2$. Then the valid $k$ should satisfy $(1 + \sqrt{2})/2 < k < (3 + 3\sqrt{2})/2$. So the valid $k$ includes only $2$. Hence, the length of intersection of the ray with pixels $(2, 2, 0)$ is $3\sqrt{2}/2 - 1$; \\
(ii) If $i = 1$, then $C_x = -1$, and the valid $j$ should satisfy $0 < j < 2$. So the valid $j$ only contains $1$. Then $C_y = 1$ for $j = 1$. Then the valid $k$ should satisfy $(1 - \sqrt{2})/2 < k < (3 + \sqrt{2})/2$. So the valid $k$ includes $0$, $1$, and $2$. The lengths of intersections of the ray with voxels $(0, 1, 1)$, $(1, 1, 1)$ and $(2, 1, 1)$ are $1 - \sqrt{2}/2$, $\sqrt{2}$, and $1 - \sqrt{2}/2$, respectively; \\
(iii) If $i = 2$, then $C_x = 1$, and then the valid $j$ should satisfy $-1 < j < 1$. Hence, the valid $j$ only contains $0$. The $C_y = 3$ for $j = 0$. Then the valid $k$ should satisfy $(1 - 3\sqrt{2})/2 < k < (3 - \sqrt{2})/2$. So the valid $k$ includes only $0$. Hence, the length of intersection of the ray with pixels $(0, 0, 2)$ is $3\sqrt{2}/2 - 1$. 

By numerical implementation, the output of the valid lengths of intersections is given in \cref{tab:test_3D_parallel}. The 2nd, 4th, 13th, 22nd and 24th voxles in \ac{1D} form are the voxels $(0, 0, 2)$, $(0, 1, 1)$, $(1, 1, 1)$, $(2, 1, 1)$ and $(2, 2, 0)$, respectively.  
  \begin{table}[htbp]
  \centering
   \begin{tabular}{cl}
   \hline
    2 & 1.12132 \\ [-0.3em]
    4 & 0.292893 \\ [-0.3em]
    13 & 1.41421 \\ [-0.3em]
    22 & 0.292893 \\ [-0.3em]
    24 & 1.12132 \\
    \hline
    \end{tabular}
    \caption{Output of the example for \ac{3D} parallel beam}
    \label{tab:test_3D_parallel}
    \end{table}
    
Evidently, the above results are completely consistent with each other if we neglect the machine error. So the error is zero.

\subsection{Test suite 4: Circular cone beam} 
\label{sec:Test_suite_4}

Since the validation of equispaced circular cone beam has no essential distinction with that of the equiangular one, here we just give an example for the equiangular case.  

As an example, the size of the image is assumed to be $4\times 4$, the given ray is assigned with parameters $D = 4$, $\phi'_1 = \pi/4$, $\alpha = \pi/12$, and $\beta = \pi/12$. Transforming these parameters into the associated ones in \ac{3D} parallel beam by \cref{eq:equiangular_cone_beam}, we have 
$\phi_1 = \pi/3$, $\phi_2 = \pi/12$, $s_1 = 4\sin(\pi/12)$, $s_2 = 4\sin(\pi/12)\cos(\pi/12)$. 

For this example, the main algorithm are almost the same as the \ac{3D} parallel beam. So we neglect the calculations by hand following the algorithm in \cref{sec:3D_cone_beam} as the test above. 
          
Using solid geometry, we can analytically calculate out the valid intersections of the ray with all the voxels. We list the valid intersected voxels and the corresponding lengths of intersections as follows  
\begin{align*}
&\text{voxel}~(1, 3, 0) : 2\bigl(1 - (4-\sqrt{2})\tan(\pi/12)\sin(5\pi/12)/\sin(\pi/3)\bigr)/\cos(\pi/12); \\[-0.5em]
&\text{voxel}~(1, 2, 0) : 2\sqrt{3}/\bigl(3\cos(\pi/12)\bigr); \\[-0.5em]
&\text{voxel}~(1, 1, 0) : 2\bigl((4-\sqrt{2})\tan(\pi/12)\sin(5\pi/12)/\sin(\pi/3) - \sqrt{3}/3\bigr)/\cos(\pi/12); \\[-0.5em]
&\text{voxel}~(0, 0, 1) : 2\sqrt{3}/\bigl(3\cos(\pi/12)\bigr); \\[-0.5em]
&\text{voxel}~(1, 1, 1) : \bigl(1 - (4\sqrt{2} - 2)\tan(\pi/12)\bigr)/\sin(\pi/12); \\[-0.5em]
&\text{voxel}~(0, 1, 1) : \bigl(\bigl((4 - \sqrt{2})\bigl(\sqrt{2} - 2\tan(\pi/12)\sin(5\pi/12)/\sin(\pi/3)\bigr) \\[-0.5em]
& ~\qquad\qquad\qquad+ 4\sqrt{3}/3\bigr)\tan(\pi/12) - 1\bigr)/\sin(\pi/12). 
\end{align*}

By numerical implementation, the output of the valid lengths of intersections is displayed in Table \cref{tab:test_circular_cone}. The 1st, 5th, 20th, 21st, 24th and 28th voxles in \ac{1D} form are the voxels $(0, 0, 1)$, $(0, 1, 1)$, $(1, 1, 0)$, $(1, 1, 1)$, $(1, 2, 0)$ and $(1, 3, 0)$, respectively.  
 \begin{table}[htbp]
  \centering
   \begin{tabular}{cl}
   \hline
    1 & 1.19543 \\ [-0.3em]
    5 & 0.712929 \\ [-0.3em]
    20 & 0.404656 \\ [-0.3em]
    21 & 0.0778492 \\ [-0.3em]
    24 & 1.19543 \\ [-0.3em]
    28 & 0.470462 \\ 
    \hline
    \end{tabular}
    \caption{Output of the example for equiangular circular cone beam}
    \label{tab:test_circular_cone}
    \end{table}
    
As we have checked up, the above results are completely consistent with each other if we neglect the machine error. So the error is zero.

\subsection{Test suite 5: Helical cone beam} 
\label{sec:Test_suite_5}

Here we just give an example for the equiangular helical cone beam. The geometry of this case is almost the same as the equiangular circular cone beam except with one more signed vertical distance. 

When we choose the vertical distance $H = 0$, and maintain those parameters in Section \cref{sec:Test_suite_4}, as expected, the output of the numerical implementation is the same as the corresponding equiangular circular cone beam. 

In contrast, we take $H = 0.5$. Transforming these parameters into the associated ones in \ac{3D} parallel beam by  \cref{eq:equiangular_cone_beam} and  \cref{eq:Helical_cone_beam_equiangular}, we have $\phi_1 = \pi/3$, $\phi_2 = \pi/12$, $s_1 = 4\sin(\pi/12)$, $s_2 = 4\cos(\pi/12)\sin(\pi/12) + \cos(\pi/12)/2$. 

Similarly, we can analytically calculate out the valid intersections of the ray with all the voxels. The results of valid intersected voxels and the corresponding lengths of intersections are listed as follows 
\begin{align*}
&\text{voxel}~(0, 2, 0) : 2\sqrt{3}/\bigl(3\cos(\pi/12)\bigr); \\[-0.5em]
&\text{voxel}~(0, 1, 0) : 2\bigl((4-\sqrt{2})\tan(\pi/12)\sin(5\pi/12)/\sin(\pi/3) - \sqrt{3}/3\bigr)/\cos(\pi/12); \\[-0.5em]
&\text{voxel}~(0, 0, 1) : 2\sqrt{3}/\bigl(3\cos(\pi/12)\bigr); \\[-0.5em]
&\text{voxel}~(0, 1, 1) : \bigl(2 - (8-2\sqrt{2})\tan(\pi/12)\sin(5\pi/12)\bigr)/\bigl(\cos(\pi/6)\cos(\pi/12)\bigr);  \\[-0.5em] 
&\text{voxel}~(0, 3, 0) : \bigl(0.5 - (4\sqrt{2} - 4)\tan(\pi/12)\bigr)/\sin(\pi/12); \\[-0.5em]
&\text{voxel}~(1, 3, 0) : 2\bigl(1 - (4-\sqrt{2})\tan(\pi/12)\sin(5\pi/12)/\sin(\pi/3)\bigr)/\cos(\pi/12)  \\[-0.5em]
&~\qquad\qquad\qquad - \bigl(0.5 - (4\sqrt{2} - 4)\tan(\pi/12)\bigr)/\sin(\pi/12). 
\end{align*}

Furthermore, the output of the numerical implementation is given in \cref{tab:test_helical_cone}. The 1st, 4th, 5th, 8th, 12th and 28th voxles in \ac{1D} form corresponds to the voxels $(0, 0, 1)$, $(0, 1, 0)$, $(0, 1, 1)$, $(0, 2, 0)$, $(0, 3, 0)$ and $(1, 3, 0)$, respectively.  
 \begin{table}[htbp]
  \centering
   \begin{tabular}{cl}
   \hline
    1 & 1.19543 \\[-0.3em]
    4 & 0.404656 \\[-0.3em]
    5 & 0.790778 \\[-0.3em]
    8 & 1.19543 \\[-0.3em]
    12 & 0.253912 \\[-0.3em]
    28 & 0.21655 \\
    \hline
    \end{tabular}
    \caption{Output of the example for equiangular helical cone beam}
    \label{tab:test_helical_cone}
    \end{table}

For this example, it is easy to validate that the analytic results and the corresponding numerical ones are also completely consistent with each other if the machine error is neglected. So the error is zero.

\section{Conclusion}
\label{sec:Conclusions}   

A fast, accurate, adaptive and parallelizable algorithm has been proposed to compute the X-ray transform of an image represented by unit (pixel/voxel) basis functions.  The algorithm is first proposed for the most basic scanning geometry (\ie \ac{2D} parallel beam), and is then extended to \ac{2D} fan beam by the associated transformation formula. The algorithm is further generalized  to \ac{3D} parallel beam, and is adapted to \ac{3D} equiangular and equispaced circular/helical cone beams by using the derived transformation formulas. The proposed algorithm is not just limited to the above scanning geometries, which is applicable to any beam that its parameters can be transformed into those of the corresponding parallel beam. Additionally, this algorithm is also applicable to the situation of the image represented by polygonal/polyhedral basis functions.  So the proposed algorithm is more flexible, \eg one can easily work with different geometries and different image basis functions. 

Moreover, the conditions for non-vanishing intersectability are derived. More importantly, for any given $i$ (or $j$), the range of valid $j$ (or $i$) is calculated by the condition in \ac{2D} situation, to obtain those pixels intersecting with the given ray non-vanishingly. In \ac{3D} situation, for any given $i$ (or $j$, $k$), the ranges of valid $j$ and $k$ (or $i$ and $k$, $i$ and $j$) are calculated by the condition, to obtain those voxels intersecting with the given ray non-vanishingly. And then the intersection lengths are computed by the obtained analytic formula for those units rather than all the individuals. The overhead that comes from evaluating the ranges is one addition and one multiplication for \ac{2D} case, and triple additions and triple multiplications for \ac{3D} case. The computational complexity is $O(1)$ from computing the analytic formula for both cases. So for any given ray, the computational cost of the proposed algorithm is $O(N)$ rather than $O(N^d)$. Here $N$ represents the size of the image being reconstructed along one axis, and $d$ denotes the dimension of image space. Hence, the algorithm becomes quite easy to be understood and implemented, and its computational cost is significantly reduced. 

To the problem itself, we further discussed the intrinsic ambiguities that perhaps happen, and have presented a solution in the algorithm. Moreover, the algorithm not only possesses the adaptability with regard to the center position, scale and size of the image, and the more general scanning geometries, but also is quite suited to parallelize with optimality. The comparison study is included in this work with respect to several existing algorithms, which demonstrates the proposed algorithm is faster and more complete to use, and is applicable to the more general case of the image represented by polygonal/polyhedral basis functions. The projection matrix can be sparsely stored and output if needed, and the adjoint of X-ray transform can be also computed by the algorithm. Hence, the algorithm can be customized freely according to the requirements of the user, and more scanning geometries can be easily added into the framework based on the proposed algorithm. 

The correctness of the proposed algorithm has been validated by several test examples in \cref{sec:Test_suite_1}--\cref{sec:Test_suite_5} for various scanning geometries. For each scanning geometry, we provided one test suite. As the pipeline, we first calculated the analytic results by solid geometry, and/or recomputed the results by hand along the proposed algorithm, which construct the standard references being used to make comparisons. And then we computed the associated numerical results by numerical implementation. As demonstrated by these compared results, the proposed algorithm can yield precise results for various scanning geometries. To our algorithm, the numerical errors do not exist other than the machine error. Hence, the error is zero. 

Based on the proposed algorithm, the future work will focus on the development of a software package for tomographic image reconstruction to various imaging modalities. And the parallelized version with CPU and GPU will also be implemented. 

\section*{Acknowledgments}
\label{sec:Acknowledgments}   

This work of Chen was partially funded by the Beijing Natural Science Foundation under grant Z180002. The work of \"Oktem was supported by the Swedish Foundation of Strategic Research grant AM13-0049.

\bibliographystyle{plain}
\bibliography{shapereferences}

\end{document}